\documentclass[11pt]{article}
\usepackage{graphicx}
\usepackage{amsfonts,bbm}
\usepackage{amsthm}
\usepackage{amssymb}
\usepackage{amsmath,amsfonts}
\usepackage{mathtools, bm}
\usepackage{xcolor}
\usepackage{dsfont}
\usepackage{comment}
\usepackage{hyperref}
\usepackage{float}
\usepackage[numbers]{natbib}
\usepackage{tikz}
\usetikzlibrary{arrows,decorations.pathreplacing, bending,hobby}
\usetikzlibrary{shapes.geometric}
\usetikzlibrary{decorations.markings}
\usepackage{subcaption}
\usepackage{cleveref}[capitalise, compress, nameinlink, noabbrev]
\usepackage[T1]{fontenc}
\usepackage[margin=1.25in]{geometry}
\usepackage{enumitem}
\usepackage{indentfirst}

\theoremstyle{plain}
\newtheorem{theorem}{Theorem}[section]
\crefname{theorem}{Theorem}{Theorems}

\newtheorem{proposition}[theorem]{Proposition}
\crefname{proposition}{Proposition}{Propositions}

\newtheorem{corollary}[theorem]{Corollary}
\crefname{corollary}{Corollary}{Corollaries}

\newtheorem{lemma}[theorem]{Lemma}
\crefname{lemma}{Lemma}{Lemmata}

\crefname{conjecture}{Conjecture}{Conjectures}

\crefname{problem}{Problem}{Problems}

\crefname{claim}{Claim}{Claims}

\crefname{observation}{Observation}{Observations}

\crefname{setup}{Setup}{Setups}

\crefname{fact}{Fact}{Facts}

\crefname{algorithm}{Algorithm}{Algorithms}

\crefname{remark}{Remark}{Remarks}

\crefname{example}{Example}{Examples}

\theoremstyle{definition}

\crefname{definition}{Definition}{Definitions}

\newtheorem{construction}[theorem]{Construction}
\crefname{construction}{Construction}{Constructions}

\crefname{question}{Question}{Questions}

\numberwithin{equation}{section}

\crefname{section}{Section}{Sections}
\crefname{appendix}{Appendix}{Appendix}

\crefname{figure}{Figure}{Figures}

\newcommand{\jan}[1]{\textcolor{orange}{\textbf{#1}}}
\newcommand{\jun}[1]{\textcolor{blue}{\textbf{#1}}}

\newcommand\ceil[1]{\left\lceil#1\right\rceil}
\newcommand\floor[1]{\left\lfloor#1\right\rfloor}
\newcommand\rc{R_{\textup{cyc}}}
\newcommand\ro{R_{\textup{ord}}}

\DeclareMathOperator{\AP}{AP}
\DeclareMathOperator{\MP}{MP}
\DeclareMathOperator{\NM}{NM}
\DeclareMathOperator{\MC}{MC}

\DeclareMathOperator{\eff}{eff}

\allowdisplaybreaks

\title{New results on cyclic and ordered Ramsey numbers} 
\author{Gaurav Kucheriya\thanks{SAFIR Institute, Sorbonne University, Abu Dhabi, United Arab Emirates. Was supported by GAČR grant 25-17377S and GAUK project 92125 during this project. Email: \url{gaurav.kucheriya@sorbonne.ae}.} \and Jan Petr\thanks{Fakultät für Informatik und Mathematik, Universität Passau, Germany. Was supported by the Deutsche Forschungsgemeinschaft (DFG, German Research Foundation) – 542321564 during this project. Email: \url{jan.petr@uni-passau.de}.} \and Jun Yan\thanks{Faculty of Mathematics, University of Waterloo, Waterloo, Canada. Was supported by ERC Advanced Grant 883810 during this project. Email: \url{j228yan@uwaterloo.ca}.}}

\begin{document}

\maketitle

\begin{abstract}
Inspired by the work of Poljak, Bašić, Damnjanović, Stevanović, and Stošić used the Kissat SAT solver in a recent paper to systematically determine the small cyclic and ordered Ramsey numbers of several families of graphs. By studying the patterns in their data, Bašić et al. formulated $15$ conjectures on cyclic and ordered Ramsey numbers. In this paper, we confirm $9$ of their conjectures, trace confirmations of $2$ more to existing literature, and, perhaps surprisingly, disprove another one. 
\end{abstract}

\section{Introduction}
An \emph{ordered graph} is a graph whose vertex set is equipped with a total order. We say that an ordered graph $G$ \emph{contains} an ordered graph $H$ if there is an injection $f:V(H) \to V(G)$ that preserves both the vertex order and adjacency, that is, for all $i < j$ in $V(H)$, we have $f(i) < f(j)$ in $V(G)$, and $f(i)f(j) \in E(G)$ whenever $ij \in E(H)$. We say that $G$ \emph{avoids} $H$ if $G$ does not contain $H$. 
Given ordered graphs $H$, $H'$, the \emph{ordered Ramsey number} $\ro(H,H')$ is defined to be the smallest natural number $N$ such that every red/blue edge-colouring of $K_N^<$, the ordered complete graph on $\{1,2,\ldots, N\}$, contains a red copy of $H$ or a blue copy of $H'$ as an ordered subgraph. When $H=H'$, we denote $\ro(H,H)$ by $\ro(H)$.

Recently, Balko, Cibulka, Král, and Kynčl~\cite{balko2013ramsey}, and independently Conlon, Fox, Lee, and Sudakov~\cite{conlon2017ordered}, systematically studied ordered Ramsey numbers. For a comprehensive survey of results in this area we refer the readers to Balko~\cite{balko2025survey}. 
In \cite{balko2013ramsey}, Balko et al. obtained the following bounds on the ordered Ramsey number of an alternating path
\[(5/2+o(1))n=5\lfloor n/2\rfloor-4\leq \ro(\AP_n)\leq 2n-3+\sqrt{2n^2-8n+11}=(2+\sqrt2+o(1))n.\] Recall that the \emph{alternating path} $\AP_n$ on $n$ vertices is the labelled path with vertex set $[n]$ such that the vertices encountered along the path are: 
$1,n,2,n-1,\ldots,n/2,n/2+1$ if $n$ is even, and $1,n,2,n-1,\ldots,\lfloor n/2 \rfloor ,\lfloor n/2\rfloor +2,\lfloor n/2\rfloor +1$ if $n$ is odd. Balko et al.~\cite[Problem 29]{balko2013ramsey} also asked whether the alternating path minimizes the ordered Ramsey number among all $n$-vertex ordered paths. This remains an open problem. 

Together with Lo and Sgueglia~\cite{kucheriya2026ordered}, we improved the upper bound on $\ro(\AP_n)$. In particular, we proved that
\[
    \ro(\AP_n)\leq2n-2+\floor{\frac{\sqrt{2(n-2)^2+(-1)^n}-1}2}=\left(2+\frac{\sqrt{2}}{2}+o(1)\right)n.
\]
In his master's thesis~\cite{poljak2023off}, Poljak listed the exact values or lower bounds for $\ro(\AP_n)$ up to $n=13$, which were obtained from computer experiments using a SAT solver.
It is surprising that our upper bound mentioned above matches all the known results reported in~\cite{poljak2023off}. However, we do not believe that our bound is sharp for large~$n$. 

Inspired by the work of Poljak~\cite{poljak2023off}, Bašić, Damnjanović, Stevanović, and Stošić~\cite{Bavsic2026SomeRO} applied the Kissat SAT solver (see~\cite{biere2024cadical,biere2024cadicalgithub}) and systematically determined new small ordered Ramsey numbers involving specific classes of graphs, including monotone paths, alternating paths, monotone cycles, stars, complete graphs, and nested matchings. Furthermore, they defined the following notion of cyclic Ramsey number, and similarly used the Kissat SAT solver to systematically determine many small cyclic Ramsey numbers.   

A \emph{cyclically ordered graph}, sometimes also called a \emph{convex geometric graph} \cite{MR2038478}, is a graph whose vertex set is equipped with a cyclic order. Equivalently, its vertices may be placed in convex position in the plane, with the cyclic order being the clockwise order of the vertices, and its edges represented by straight-line segments. Two linear representations of a cyclic order are regarded as equivalent if one is obtained from the other by a cyclic shift. Under the name of convex geometric graphs, Tur\'an and Ramsey type problems for cyclically ordered graphs had been studied before, for example, in~\cite{MR2038478,MR3697924}.

We say that a cyclically ordered graph $G$ \emph{contains} a cyclically ordered graph $H$ if there is an injection $f\colon V(H)\to V(G)$ that preserves both adjacency and the cyclic order. More explicitly, if the vertices of $H$ and $G$ are placed in convex position, and the vertices in $H$ appear as $v_1,v_2,\ldots,v_h$ in cyclic order, then the vertices $f(v_1),\ldots,f(v_h)$ also appear in cyclic order in $G$, and $f(v_i)f(v_j)\in E(G)$ whenever $v_iv_j\in E(H)$. Thus, cyclic shifts are permitted, whereas reversing the cyclic order is not. If no such injection exists, we say that $G$ \emph{avoids} $H$.

Given cyclically ordered graphs $H$ and $H'$, the \emph{cyclic Ramsey number} $R_{\mathrm{cyc}}(H,H')$ is the smallest natural number $N$ such that every red/blue edge-colouring of $K_N^{\mathrm{cyc}}$, the complete graph on $N$ cyclically ordered vertices, contains a red cyclically ordered copy of $H$ or a blue cyclically ordered copy of $H'$. 

For both the ordered and cyclic Ramsey numbers, Bašić, Damnjanović, Stevanović, and Stošić~\cite{Bavsic2026SomeRO} also used their data to formulate 15 conjectures regarding their exact values. In this paper, we confirm 9 of their conjectures, trace confirmations of $2$ more to existing literature, and, perhaps surprisingly, disprove another conjecture of theirs. Below is a table summarizing the current status of their conjectures. The graphs involved in these conjectures are defined as follows.

$\MP_n$ and $\MC_n$ denote the \textit{monotone path} and the \textit{monotone cycle} with vertex set $[n]$, respectively, where in both cases the vertices encountered along the path or cycle are $1,2,\ldots,n$ in order. $\AP^r_n$ denotes the \textit{reverse alternating path} with vertex set $[n]$, which is obtained from the alternating path $\AP_n$ defined above by reversing the vertex order. $K_n$ denotes the complete graph on $n$ vertices, while $S_n$ denotes the \textit{start-central star} with vertex set $[n]$, which is the star whose centre is the vertex 1. Finally, $\NM_{2n}$ denotes the \textit{nested matching} with vertex set $[2n]$, whose set of $n$ edges is $\{\{i,2n+1-i\}:i\in[n]\}$.

\begin{table}[h!]
\begin{center}\renewcommand{\arraystretch}{1.3}
\resizebox{\textwidth}{!}{%
\begin{tabular}{|c c c c c c c c|}
\hline
Conjecture & Type & $G_1$ & $G_2$ & Conjectured value & Conditions & Status & Reference\\
\hline\hline
$4.7$ & $\rc$ & $\MP_k$ & $\MP_\ell$ & $1+(k-1)(\ell-2)$ & $\ell\geq k \geq 3$ & True & \cite[Theorem 1]{MR3697924}\\ 
\hline
$4.8$ & $\rc$ & $\AP_k$ & $\AP_\ell$ & $k+\ell-2-(k\ell \mod 2)$ & $k\geq 2, \ell\geq 2$ & True & Theorem~\ref{thm:cycapap}\\
\hline
$4.9$ & $\rc$ & $\AP_k$ & $\MP_\ell$ & $1+(k-1)(\ell-2)$ & $k\geq 3, \ell\geq 4$ & True & Theorem~\ref{thm:cycapmp}\\
\hline
$4.10$ & $\rc$ & $\AP_k$ & $\AP^r_\ell$ & $k+\ell-2-(k\ell \mod 2)$ & $k \geq 2, \ell \geq 2$ & True & Theorem~\ref{thm:cycapap}\\
\hline
$4.13$ & $\ro$ & $\MC_k$ & $\AP_\ell$ & $\left\lceil (k-\frac{1}{2})(\ell-1)\right\rceil$ & $k\geq 3, \ell\geq 2$ & Open & \\
\hline
$4.14$ & $\rc$ & $\MC_k$ & $\AP_\ell$ & $1+(k-1)(\ell-1)$ & $k\geq2,\ell\geq 1$ & True & Theorem~\ref{thm:kmcap}\\
\hline
$4.16$ & $\rc$ & $S_k$ & $\MP_\ell$ & $1+(k-1)(\ell-2)$ & $k\geq 3, \ell\geq 4$ & True & \cite[Theorem 4]{MR3697924}\\
\hline
$4.17$ & $\rc$ & $S_k$ & $\AP_\ell$ & $k+\ell-2-(k\ell \mod 2)$ & $k,\ell\geq 2$ & True & Theorem~\ref{thm:cycapap}\\
\hline
$4.22$ & $\ro$ & $K_k$ & $\AP_\ell$ & $\left\lceil \frac{3k\ell-5k}{2} \right\rceil - 2\ell+5$ & $k\geq 2, \ell\geq 3$ & Open & \\
\hline
$4.23$ & $\rc$ & $K_k$ & $\AP_\ell$ & $1+(k-1)(\ell-1)$ & $k\geq1, \ell\geq1$ & True & Theorem~\ref{thm:kmcap}\\
\hline
$4.30$ & $\rc$ & $\NM_{2k}$ & $\NM_{2\ell}$ & $2k+2\ell-3$ & even $k,\ell\geq2$ & True & Theorem~\ref{thm:cycnmnm}\\
\hline
$4.31$ & $\ro$ & $\NM_{2k}$ & $\MP_\ell$ & $1+(2k-1)(\ell-1)$ & $k,\ell\geq 1$ & True & Theorem~\ref{thm:ordnmmp}\\
\hline
$4.32$ & $\rc$ & $\NM_{2k}$ & $\MP_\ell$ & $k\ell-k+1-(\ell-2)\mathbbm{1}\{2\mid k\}$ & $k\geq 1, \ell \geq 3$ & LB & Proposition~\ref{prop:lower}\\
\hline
$4.33$ & $\rc$ & $\NM_{2k}$ & $S_\ell$ & $2k+\ell-2-\mathbbm{1}\{2\mid k, 2 \nmid \ell\}$ & $k \geq 1, \ell\geq 2$ & True & Theorem~\ref{thm:cycnms}\\
\hline
$4.36$ & $\rc$ & $\NM_{2k}$ & $K_\ell$ & $k\ell$ & $k\geq 1, \ell\geq 2$ & False & Theorem~\ref{thm:disproof}\\
\hline
\end{tabular}
}
\caption{Conjectures from \cite{Bavsic2026SomeRO} and their status. $\mathbbm{1}\{\cdot\}$ is the indicator function that takes value 1 when the statement within holds, and 0 otherwise. Status LB of Conjecture~4.32 denotes that so far it is known only that the conjectured value is a lower bound on $\rc(\NM_{2k},\MP_{\ell})$.}\label{table}
\end{center}
\end{table}

We remark that from our results together with the data and conjectures of Bašić et al.~\cite{Bavsic2026SomeRO}, cyclic Ramsey numbers generally seem to admit simpler formulas compared to their ordered counterparts. Intuitively, this might be explained by the additional rotational symmetry. 

It turns out that Conjecture~4.7 and Conjecture~4.16 have already been proved back in 2017 by Jayawardene, Rousseau, and Harborth in Theorem 1 and Theorem 4 of~\cite{MR3697924}, respectively. The rest of the paper presents our remaining results. In Section~\ref{sec:cp}, we confirm Conjectures~4.8, 4.9, 4.10, 4.14, 4.17, and 4.23, which will follow from more general results on a family of trees called caterpillars. In Section~\ref{sec:nm}, which is devoted to cyclic and ordered Ramsey numbers involving nested matchings, we confirm Conjectures~4.30, 4.31, and 4.33, as well as establishing the lower bound part of Conjecture~4.32. Finally, in Section~\ref{sec:disprove}, we disprove Conjecture~4.36 by generalizing a construction of Cs\'oka, Bl\'azsik, Kir\'aly, and Lenger~{\cite{MR4442957}} in the context of necklace folding problems.

\section{Exact cyclic Ramsey numbers involving caterpillars}\label{sec:cp}
In this section, we confirm Conjectures 4.8, 4.9, 4.10, 4.14, 4.17, and 4.23 from \cite{Bavsic2026SomeRO}. All of these conjectures are about cyclic Ramsey numbers involving cyclically ordered caterpillars. A \textit{caterpillar} is a tree containing a central path, such that every vertex in the tree is either on the central path, or adjacent to a vertex on the central path. In particular, paths and stars are caterpillars.  

We will use the following Tur\'an-type result due to Perles recorded on page 292 of~\cite{moser1993recent}, reformulated using our terminology.\footnote{Note that in \cite{moser1993recent}, caterpillars are defined differently, namely as \emph{forests} that do not contain three paths with three vertices each that share one common endpoint. 
For trees, these two definitions are equivalent.} For the proof of this lemma (in fact, of a more general result), the reader is referred to Theorem 7 in Brass, Károlyi and Valtr~\cite{MR2038478}.

\begin{lemma}[Perles, 1990]\label{lemma:cycbipturan}
Let $F$ be a cyclically ordered caterpillar with $k\geq 2$ vertices. Then, for every $N>k$, the maximum number of edges an $N$-vertex cyclically ordered graph with no cyclic copy of $F$ can have is $\lfloor(k-2)N/2\rfloor$.
\end{lemma}

We also use the following two notions of distances in cyclically ordered graphs. Let $G$ be the cyclically ordered complete graph on $N$ vertices, arranged on a circle, so that vertices $1,2,\ldots,N$ appear in clockwise order. 

The \textit{distance} between two vertices $i,j\in[N]$ is defined to be $d(i,j)=\min\{|i-j|,N-|i-j|\}$. For a cyclically ordered subgraph $H$ of $G$, the \textit{distance set} of $H$, denoted by $D(H)$, is defined to be the set $\{d(i,j):ij\in E(H)\}$.

The \textit{clockwise distance} between two vertices $i,j\in[N]$, denoted by $d_c(i,j)$, is 1 plus the number of vertices on the clockwise arc $ij$, excluding $i$ and $j$. Note that we always have $d(i,j)=\min\{d_c(i,j),N-d_c(i,j)\}$.
\subsection{Alternating paths versus caterpillars}
In this subsection, we confirm Conjectures 4.8, 4.10, and 4.17 from \cite{Bavsic2026SomeRO}. That is, we obtain the exact cyclic Ramsey numbers of $\AP_{\ell}$ versus $\AP_k$, $\AP_k^r$, and $S_k$. In all three cases, the upper bound follows from Lemma~\ref{lemma:cycbipturan}, while the lower bound constructions are given by Proposition~\ref{prop:aplower}. To prove Proposition~\ref{prop:aplower}, we need the following lemma that controls the distance set sizes.
\begin{lemma}\label{lemma:distance}
Let $N\geq k\geq 2$ be integers. Then, for any cyclic copy $P$ of $\AP_k$, $\AP^r_k$, or $S_k$ in the cyclically ordered complete graph on $N$ vertices, each distance in $D(P)$ can be achieved by at most two edges in $P$. Moreover, if $N$ is even, then the distance $N/2$ is achieved by at most one edge in $P$. 

Therefore, $D(P)$ has size at least $\ceil{(k-1)/2}$. Furthermore, if $N$ is even and $D(P)$ contains $N/2$, then $D(P)$ has size at least $\ceil{k/2}$.
\end{lemma}
\begin{proof}
First let $P$ be any cyclic copy of $\AP_k$. Let the vertices in $P$ be $v_1,v_2,\ldots,v_k$, in order of their appearances on $P$, so that $v_1$ corresponds to the vertex 1 in $\AP_k$. From the definition of $\AP_k$, all of $v_3,\ldots,v_k$ lie on the clockwise arc $v_2v_1$, and the $k-1$ clockwise distances $d_c(v_2,v_1),d_c(v_2,v_3),d_c(v_4,v_3),d_c(v_4,v_5),\ldots$ are strictly decreasing. Along with $d(v_i,v_{i+1})=\min\{d_c(v_i,v_{i+1}),N-d_c(v_i,v_{i+1})\}$, it follows that for every $1\leq d\leq\floor{N/2}$, there are at most two indices $i\in[k-1]$ satisfying $d(v_i,v_{i+1})=d$. Furthermore, if $N$ is even, there is at most one index $i\in[k-1]$ satisfying $d(v_i,v_{i+1})=N/2$. Therefore, $|D(P)|\geq\ceil{(k-1)/2}$, and if $N$ is even and $N/2\in D(P)$, then $|D(P)|\geq1+\ceil{(k-2)/2}=\ceil{k/2}$. 

The case when $P$ is a cyclic copy of $\AP^r_k$ follows by reflection. Now, if $P$ is any cyclic copy of $S_k$, let $v_1$ be the centre of the star, and let $v_2,\ldots,v_k$ be the leaves in clockwise order. Then, the $k-1$ clockwise distances $d_c(v_1,v_2),d_c(v_1,v_3),\ldots,d_c(v_1,v_k)$ are strictly increasing, so we can use the same argument as above to finish the proof.
\end{proof}

Now we can give the lower bound constructions. They differ slightly depending on parities, but in all cases, they have the property that the colour of each edge $ij$ depends only on the distance $d(i,j)$. 
\begin{proposition}\label{prop:aplower}
For all integers $k,l\geq2$ and for each $T_k\in\{\AP_k,\AP^r_k,S_k\}$, we have $$\rc(T_k,\AP_\ell)\geq k+\ell-2-(k\ell \bmod 2).$$
\end{proposition}
\begin{proof}
Let $N=k+\ell-3-(k\ell \bmod 2)$, it suffices to find a red/blue colouring of the edges of the cyclically ordered complete graph $G$ on $N$ vertices so that there is no red copy of $T_k$ and no blue copy of $\AP_\ell$.
Note that the distance between any two distinct vertices of $G$ lies in the set $D=\{1,2,\ldots,\floor{N/2}\}$. 

If both $k$ and $\ell$ are even, then $N=k+\ell-3$ is odd, and we can partition $D$ as $K\cup L$, so that $|K|=(k-2)/2$ and $|L|=(\ell-2)/2$. For every edge $uv$ in $G$, if $d(u,v)$ lies in $K$, then colour it red, otherwise colour it blue. Then, by Lemma~\ref{lemma:distance}, there is no red copy of $T_k$ and no blue copy of $\AP_\ell$. 

If both $k$ and $\ell$ are odd, then $N=k+\ell-4$ is even, and we can partition $D\setminus\{N/2\}$ as $K\cup L$, so that $|K|=(k-3)/2$ and $|L|=(\ell-3)/2$. For every edge $uv$ in $G$, colour it red if $d(u,v)$ lies in $K\cup\{N/2\}$, and colour it blue if $d(u,v)$ lies in $L$. Then, by Lemma~\ref{lemma:distance}, there is no red copy of $T_k$ and no blue copy of $\AP_\ell$. 

If exactly one of $k$ and $\ell$ is odd, say $k$ is odd, then $N=k+\ell-3$ is even, and we can partition $D\setminus\{N/2\}$ as $K\cup L$, so that $|K|=(k-3)/2$ and $|L|=(\ell-2)/2$. For every edge $uv$ in $G$, colour it red if $d(u,v)$ lies in $K\cup\{N/2\}$, and colour it blue if $d(u,v)$ lies in $L$. Then, by Lemma~\ref{lemma:distance}, there is no red copy of $T_k$ and no blue copy of $\AP_\ell$. 
\end{proof}

Since paths and stars are caterpillars, we can now apply Lemma~\ref{lemma:cycbipturan} to confirm Conjectures~4.8, 4.10, and 4.17 from \cite{Bavsic2026SomeRO} as follows.
\begin{theorem}\label{thm:cycapap}
For all integers $k,l\geq2$, we have
$$\rc(\AP_k,\AP_\ell)=\rc(\AP_k,\AP^r_\ell)=\rc(S_k,\AP_\ell)=k+\ell-2-(k\ell \bmod 2).$$
\end{theorem}
\begin{proof}
The lower bound follows from Proposition~\ref{prop:aplower}. For the upper bound, let $G$ be a red/blue coloured cyclically ordered complete graph on $N=k+\ell-2-(k\ell \bmod 2)$ vertices. Let $T_k\in\{\AP_k,S_k\}$, and let $T_\ell'\in\{\AP_\ell,\AP^r_\ell\}$. By Lemma~\ref{lemma:cycbipturan}, if $G$ contains no red cyclic copy of $T_k$ and no blue cyclic copy of $T_\ell'$, then the number of red edges in $G$ is at most $\floor{(k-2)N/2}$, and the number of blue edges in $G$ is at most $\floor{(\ell-2)N/2}$. However, it is straightforward to check that $\binom{N}{2}>\floor{(k-2)N/2}+\floor{(\ell-2)N/2}$, a contradiction.
\end{proof}

\subsection{Alternating paths versus complete graphs and monotone cycles}
In this subsection, we confirm Conjectures~4.14 and 4.23 from \cite{Bavsic2026SomeRO}. That is, we prove that both $\rc(K_k,\AP_\ell)$ and $\rc(\MC_k,\AP_\ell)$ are equal to $1+(k-1)(\ell-1)$. Note that $\MC_k$ is a subgraph of $K_k$, so we have $\rc(K_k,\AP_\ell)\geq\rc(\MC_k,\AP_\ell)$. Thus, it suffices to prove that $\rc(\MC_k,\AP_\ell)\geq 1+(k-1)(\ell-1)$ and $\rc(K_k,\AP_\ell)\leq1+(k-1)(\ell-1)$. 

The lower bound was already obtained in~\cite{Bavsic2026SomeRO} but without a detailed proof since the same construction idea was already used in that paper to prove another result. We thus include a proof here for completeness. The construction partitions the vertex set into consecutive blocks, colours edges within blocks with one colour, and edges across blocks with another. 
\begin{lemma}[{\cite[Theorem 4.6]{Bavsic2026SomeRO}}]\label{lemma:4.6}
Let $H$ be a cyclically ordered connected graph with $\ell\geq 1$ vertices. Then, for any $k\geq 2$, we have
$$\rc(\MC_k,H)\geq 1+(k-1)(\ell-1).$$
\end{lemma}
\begin{proof}
Let $N=(k-1)(\ell-1)$, and consider the red/blue coloured complete graph $K_N$ defined as follows. For each $i\in[k-1]$, let $B_i$ be the set of $\ell-1$ vertices $\{(i-1)(\ell-1)+1,(i-1)(\ell-1)+2,\ldots,i(\ell-1)\}$, and colour all edges within $B_i$ blue. For any $1\leq i<j\leq k-1$, colour all edges between $B_i$ and $B_j$ red. Then, there is no blue cyclic copy of $H$, as each blue connected component has size $\ell-1<\ell$. If there is a red cyclic copy of $\MC_k$, say with vertices $v_1,\ldots,v_k$ in cyclic order, then by pigeonhole, there exists $i\in[k-1]$ and $j\in[k]$, such that $v_j$ and $v_{j+1}$ are both in $B_i$, where we interpret $v_{k+1}$ as $v_1$. This is a contradiction as $v_jv_{j+1}$ is blue from the definition. Therefore, $\rc(\MC_k,H)\geq N+1=1+(k-1)(\ell-1)$.
\end{proof}

For the upper bound, we again use a Turán-type counting argument.
\begin{lemma}\label{lemma:completeap}
For all $k\geq2$ and $\ell \geq 1$, we have
$$\rc(K_k,\AP_{\ell})\leq 1+(k-1)(\ell-1).$$ 
\end{lemma}
\begin{proof}
Let $G$ be any red/blue coloured cyclically ordered complete graph on $N=1+(k-1)(\ell-1)$ vertices. If $G$ contains no red cyclic copy of $K_k$, then $G$ also contains no ordinary unordered copy of $K_k$, so by Tur\'an's Theorem, the number of red edges in $G$ is at most $\binom{k-2}2(\ell-1)^2+(k-2)\ell(\ell-1)$. If $G$ contains no blue cyclic copy of $\AP_\ell$, then by Lemma~\ref{lemma:cycbipturan}, the number of blue edges in $G$ is at most $\floor{(\ell-2)N/2}$. However, 
\[\binom{N}2>\floor{(\ell-2)N/2}+\binom{k-2}2(\ell-1)^2+(k-2)\ell(\ell-1),\] so $G$ must contain either a red cyclic copy of $K_k$, or a blue cyclic copy of $\AP_\ell$.
\end{proof}

Combining the two previous lemmas, we can now easily confirm Conjectures 4.14 and 4.23 from \cite{Bavsic2026SomeRO} as follows.
\begin{theorem}\label{thm:kmcap}
Let $k,\ell\geq 1$ and $k'\geq 2$. We have 
$$\rc(\MC_{k'},\AP_\ell)=1+(k'-1)(\ell-1) \quad \text{and} \quad \rc(K_k,\AP_\ell)=1+(k-1)(\ell-1).$$ 
\end{theorem}
\begin{proof}
The $k=1$ case is trivial. For $k\geq2$, since $\MC_k$ is a subgraph of $K_k$, we have $\rc(K_k,\AP_\ell)\geq\rc(\MC_k,\AP_\ell)$. By Lemma~\ref{lemma:4.6} and Lemma~\ref{lemma:completeap}, 
\[1+(k-1)(\ell-1)\leq\rc(\MC_k,\AP_\ell)\leq\rc(K_k,\AP_\ell)\leq1+(k-1)(\ell-1),\]
and so equalities hold throughout. 
\end{proof}

\subsection{Caterpillars versus monotone paths}
The goal of this subsection is to confirm Conjecture~4.9 from \cite{Bavsic2026SomeRO} by showing $\rc(\AP_k,\MP_\ell)=1+(k-1)(\ell-2)$. 
In that paper, the following general lower bound for cyclic Ramsey numbers of monotone paths versus arbitrary cyclically ordered connected graphs was established, using the same block-type construction as in Lemma~\ref{lemma:4.6}.

\begin{lemma}[{\cite[Theorem 4.5]{Bavsic2026SomeRO}}]\label{lemma:monogen}
Let $k,\ell$ be positive integers. Let $H$ be a $k$-vertex cyclically ordered connected graph. Then,
\[\rc(H,\MP_\ell)\geq1+(k-1)(\ell-2).\]
\end{lemma}

To confirm Conjecture 4.9, we also need the following result of Jayawardene, Rousseau, and Harborth, which actually confirms Conjecture~4.16 as mentioned in the introduction.
\begin{theorem}[{\cite[Theorem 4]{MR3697924}}]\label{thm:mindegmp}
For all integers $k\geq3$ and $\ell\geq4$, we have $$\rc(S_k,\MP_\ell)=1+(k-1)(\ell-2).$$
\end{theorem}

We now prove a more general result that when $F$ is a cyclically ordered caterpillar, the lower bound in Lemma~\ref{lemma:monogen} gives the correct value. In particular, this confirms Conjecture~4.9 from \cite{Bavsic2026SomeRO} since alternating paths are cyclically ordered caterpillars. In fact, by substituting Lemma~\ref{lemma:cycbipturan} with Theorem 7 from \cite{MR2038478} in the following proof, it is possible to show that Theorem~\ref{thm:cycapmp} holds for a more general class of cyclically ordered graphs, which we do not introduce here for simplicity of presentation.

\begin{theorem}\label{thm:cycapmp}
Let $k\geq3$ and $\ell\geq4$ be integers and let $F$ be a cyclically ordered caterpillar on $k$ vertices. Then
$$\rc(F,\MP_\ell)=1+(k-1)(\ell-2).$$
\end{theorem}
\begin{proof}
The lower bound follows from Lemma~\ref{lemma:monogen}, so it suffices to prove the upper bound. 

Let $G$ be a red/blue coloured cyclically ordered complete graph on $N=1+(k-1)(\ell-2)$ vertices. Suppose for a contradiction that $G$ contains no red cyclic copy of $F$, and no blue cyclic copy of $\MP_\ell$. Then, by Lemma~\ref{lemma:cycbipturan}, the number of red edges in $G$ is at most $\floor{(k-2)N/2}$.

By Theorem~\ref{thm:mindegmp}, if $j\geq2$, $\ell\geq 4$, and $m\geq1+j(\ell-2)$, then any red/blue coloured cyclically ordered complete graph on $m$ vertices must either contain a red vertex with degree at least $j$, and so a red cyclic copy of $S_{j+1}$, or a blue cyclic copy of $\MP_\ell$. Note that the same is trivially true for $j=1$ and $m\geq\ell$. It follows that for any nonempty induced subgraph $H$ of $G$ with $|H|\geq\ell$, there must be a vertex in $H$ with red degree at least $\floor{(|H|-1)/(\ell-2)}$ in $H$, since $G$ and thus $H$ contain no blue cyclic copy of $\MP_\ell$.

Consider now the following algorithm. We initialize by setting $G_0=G$. For every $0\leq i\leq N-\ell=(k-2)(\ell-2)-1$, suppose we have obtained an induced subgraph $G_i$ of $G$ with $N-i$ vertices. Since $N-i\geq\ell$, from above we have that $G_i$ contains a vertex $v_i$ with red degree at least $\floor{(|G_i|-1)/(\ell-2)}$ in $G_i$. Let $G_{i+1}=G_i-v_i$, and continue to the next step of the algorithm. Note that in step $i$, the deletion of $v_i$ removes at least $\floor{(|G_i|-1)/(\ell-2)}=\floor{((k-1)(\ell-2)-i)/(\ell-2)}$ red edges from $G$, and these red edges are distinct across steps. Therefore, the number of red edges in $G$ is at least
\begin{align*}
\sum_{i=0}^{(k-2)(\ell-2)-1}\floor{\frac{(k-1)(\ell-2)-i}{\ell-2}}&=(k-1)+(\ell-2)\sum_{j=2}^{k-2}j+(\ell-3)\cdot1\\
&=(k-2)(N+1)/2>\floor{(k-2)N/2},
\end{align*}
a contradiction. 
\end{proof}

\section{Ramsey numbers involving nested matchings}\label{sec:nm}
In this section, we confirm Conjectures~4.30, 4.31, and 4.33, as well as establishing the lower bound part of Conjecture~4.32. All of these results involve nested matchings. First, we prove a result on the edge distances of any cyclic copy of $\NM_{2k}$.
\begin{lemma}\label{lemma:nmdistance}
Let $G$ be an $N$-vertex cyclically ordered graph with vertices $1,2,\ldots,N$ arranged in clockwise order on a circle. Let $M$ be a cyclic copy of $\NM_{2k}$ with vertices $a_1,\ldots,a_k,b_k,\ldots,b_1$ in cyclic order, such that the $k$ edges in $M$ are $a_1b_1,\ldots,a_kb_k$. For every $i\in[k]$, let $d_i=d(a_i,b_i)$. Then, the following hold. 
\begin{enumerate}[label=\alph*)]
    \item\label{item:1} The sequence $d_1,\ldots,d_k$ is unimodal, so for each $1\leq d\leq\floor{N/2}$, there are at most two indices $i\in[k]$ such that $d_i=d$.
    \item\label{item:2} Let $d_{\max}=\max\{d_i:i\in[k]\}$. If $d_i<d_{i+1}\leq d_{\max}$ or $d_{\max}\geq d_i>d_{i+1}$, then $|d_i-d_{i+1}|\geq2$.
    \item\label{item:3} $d(a_1,a_k),d(b_1,b_k)\geq k-1$.
\end{enumerate}
\end{lemma}
\begin{proof}
For each $i\in[k]$, let $d_i'=d_c(a_i,b_i)$ be the clockwise distance from $a_i$ to $b_i$, and recall that $d_i=\min\{d_i',N-d_i'\}$. Note that for each $i\in[k-1]$, $d_{i+1}'=d_i'-d_c(a_i,a_{i+1})-d_c(b_{i+1},b_i)\leq d_i'-2$. Thus, $d_1',\ldots,d_k'$ is strictly decreasing, and so $d_1,\ldots,d_k$ is unimodal. In particular, each distance can appear at most twice in the sequence, proving~\ref{item:1}. Moreover, if $d_i<d_{i+1}\leq d_{\max}$ or $d_{\max}\geq d_i>d_{i+1}$, then either $d_i=d_i'$ and $d_{i+1}=d_{i+1}'$, or $d_i=N-d_i'$ and $d_{i+1}=N-d_{i+1}'$. In both cases, $|d_i-d_{i+1}|=d_i'-d_{i+1}'\geq2$, proving~\ref{item:2}.

For~\ref{item:3}, note that from definition, there are at least $k-2$ vertices on the clockwise arc strictly between $a_1$ and $a_k$, namely $a_2,\ldots,a_{k-1}$. Thus, $d_c(a_1,a_k)$ and similarly $d_c(b_k,b_1)$ are both at least $k-1$. Since $d_c(a_k,a_1)\geq k+1$ as the clockwise arc contains all of $b_1,\ldots,b_k$, we have $d(a_1,a_k)\geq k-1$, and similarly $d(b_1,b_k)\geq k-1$. 
\end{proof}

\subsection{Cyclic nested matchings versus stars}
In this subsection, we confirm Conjecture~4.33 from \cite{Bavsic2026SomeRO}. We begin by determining the cyclic Tur\'an number of $\NM_{2k}$, though only the upper bound will be needed.
\begin{proposition}\label{prop:nmcycturan}
Let $N\geq 2k\geq2$. The maximum number of edges in an $N$-vertex cyclically ordered graph $G$ without a cyclic copy of $\NM_{2k}$ is $(k-1)N$.
\end{proposition}
\begin{proof}
Let $G$ be an $N$-vertex cyclically ordered graph without a cyclic copy of $\NM_{2k}$. As usual, we view the $N$ vertices as $1,2,\ldots,N$ arranged clockwise on a circle, and view the vertex labels modulo $N$. If $N$ is odd, consider the edge sets $E_i=\{\{i-j,i+j\}:1\leq j\leq (N-1)/2\}$ for every $i\in[N]$, noting that they form a partition of the edge set of the $N$-vertex cyclically ordered complete graph. Note that each $E_i$ itself is a cyclic nested matching, so if $G$ contains no cyclic copy of $\NM_{2k}$, then $G$ can contain at most $k-1$ edges from each $E_i$, and hence $e(G)\leq(k-1)N$. If $N$ is even, consider the edge sets $E_i=\{\{i-j,i+j\}:1\leq j\leq(N-2)/2\}$ and $F_i=\{\{i-j,i+j-1\}:1\leq j\leq N/2\}$ for every $i\in[N/2]$. Again, these $N$ sets partition the edge set of the $N$-vertex cyclically ordered complete graph, and each of these $N$ sets is itself a cyclic nested matching. Therefore, $G$ can contain at most $k-1$ edges from each $E_i$ and each $F_i$, so $e(G)\leq(k-1)N$. This proves the upper bound. For the lower bound, we use the following constructions depending on the parities of $k$ and $N$.

If $k$ is odd, let $G$ be the $N$-vertex ordered graph containing all edges $ij$ with $d(i,j)\leq k-1$, so that $G$ is $(2k-2)$-regular and contains exactly $(k-1)N$ edges. By Lemma~\ref{lemma:nmdistance}\ref{item:1} and~\ref{item:2}, if $G$ contains a cyclic copy of $\NM_{2k}$, then there must be an edge $ij$ in this copy with $d(i,j)\geq1+2\cdot(k-1)/2=k>k-1$, a contradiction. 

If $k$ is even and $N=2n+1$ is odd, let $G$ be the $N$-vertex ordered graph such that $ij$ is an edge if and only if $n-k+2\leq d(i,j)\leq n$, so that there are $(k-1)N$ edges. Then, there is no cyclic copy of $\NM_{2k}$. Indeed, if such a copy contains no edge with distance $n$, then by Lemma~\ref{lemma:nmdistance}\ref{item:1} and~\ref{item:2}, there must be an edge $ij$ with $d(i,j)\leq n-1-2(k/2-1)=n-k+1$, a contradiction. If such a copy contains a distance $n$ edge, then it must contain exactly one as any two distance $n$ edges intersect. Then, by Lemma~\ref{lemma:nmdistance}\ref{item:1} and~\ref{item:2} again, there must be an edge $ij$ with $d(i,j)\leq n-2(k/2)=n-k$, a contradiction. 

If $k$ is even and $N=2n$ is also even, let $G$ be the $N$-vertex ordered graph such that $ij$ is an edge if and only if $n-k+2\leq d(i,j)\leq n$ or $d(i,j)=n-k+1$ and $i\in[n]$, so that there are $N/2+(k-2)N+N/2=(k-1)N$ edges. We claim that there is no cyclic copy of $\NM_{2k}$ in $G$. Like above, any such copy contains at most one distance $n$ edge. Any such copy must also contain at most one distance $n-k+1$ edge. Indeed, if it contains two, then they must be the two outermost edges in the copy of $\NM_{2k}$. However, even in the best case where these two edges are $\{1,n-k+2\}$ and $\{n,2n-k+1\}$, we have $d(n-k+2,n)=k-2<k-1$, contradicting Lemma~\ref{lemma:nmdistance}\ref{item:3}. Then, by Lemma~\ref{lemma:nmdistance}\ref{item:1} and~\ref{item:2}, the number of edges in this copy is at most $(n-(n-k))/2+((n-1)-(n-k+1))/2=k-1$, a contradiction. 
\end{proof}

We can now confirm~{\cite[Conjecture 4.33]{Bavsic2026SomeRO}}. 
The proof of the $2k+\ell-2$ case is by a simple edge counting argument. The $2k+\ell-3$ case is a bit more complicated as we need to show that the extremal configuration cannot happen. The idea is similar to the one used to prove that the Cayley table (Latin square) of $\mathbb{Z}_{2n}$ has no transversal.

\begin{theorem}\label{thm:cycnms}
For all integers $k\geq1$ and $\ell\geq2$, we have
\[\rc(\NM_{2k},S_\ell)=\begin{cases}2k+\ell-3, &\text{if }2\mid k\text{ and }\ell\text{ is odd},\\2k+\ell-2, &\text{otherwise}.
\end{cases}\]
\end{theorem}
\begin{proof}
The lower bounds are given by~{\cite[Proposition 4.34, Proposition 4.35]{Bavsic2026SomeRO}}, so we only need to prove the upper bounds. 

First suppose either $2\nmid k$ or $\ell$ is even. Let $N=2k+\ell-2$, and let $G$ be a red/blue coloured $N$-vertex cyclically ordered complete graph. If $G$ contains no red cyclic copy of $\NM_{2k}$ and no blue cyclic copy of $S_\ell$, then the number of red edges in $G$ is at most $(k-1)N$ by Proposition~\ref{prop:nmcycturan}, and every vertex in $G$ has blue degree at most $\ell-2$, so the number of blue edges in $G$ is at most $\floor{(\ell-2)N/2}$. However, $\binom{N}{2}>(k-1)N+\floor{(\ell-2)N/2}$, so $G$ must contain either a red copy of $\NM_{2k}$ or a blue copy of $S_\ell$.  

Now suppose $2\mid k$ and $\ell$ is odd. Let $N=2k+\ell-3$, and let $G$ be any red/blue coloured $N$-vertex cyclically ordered complete graph. In this case, $N$ is even, and $\binom{N}{2}=(k-1)N+\floor{(\ell-2)N/2}$. Suppose for a contradiction that $G$ contains no red cyclic copy of $\NM_{2k}$, and no blue cyclic copy of $S_\ell$. Then, since equality holds, each vertex in $G$ must have blue degree exactly $\ell-2$, and thus red degree exactly $2k-2$. Furthermore, for every $i\in[N/2]$, each set $E_i$ and $F_i$ as in the proof of Proposition~\ref{prop:nmcycturan} must contain exactly $k-1$ red edges. 

Let $N=2n$. To show that this is impossible, we relabel the $N$ vertices in $G$ as $1,\overline{1},2,\overline{2},\ldots,n,\overline{n}$ in clockwise order. Then, after updating the vertex labels, the $N$ sets $E_1,\ldots,E_{N/2},F_1,\ldots,F_{N/2}$ becomes (after reindexing) $E_1',\ldots,E_{n}',F_1',\ldots,F_{n}'$, where $E_i'=\{a\overline{b}:a+b=i\}$ and $F_i'=\{ab:a+b=i+1,a\not=b\}\cup\{\overline{a}\overline{b}:a+b=i,a\not=b\}$. Here and from now on, addition is modulo $n$. These $2n$ new sets still partition the edge set of the cyclically ordered complete graph on $N$ vertices, and the red subgraph of $G$ still contains exactly $k-1$ edges from each $E_i'$ and each $F_i'$.

Let $X,Y,Z$ be the sets of red edges in $G$ of the form $ab,\overline{a}\overline{b},a\overline{b}$, respectively. Then, $|X|+|Y|+|Z|=(k-1)N$ since the red subgraph of $G$ is $(2k-2)$-regular. By summing the degrees of all vertices of the form $a$, we have $2|X|+|Z|=(2k-2)n$. Since all edges in $E_1',\ldots,E_n'$ are of the form $a\overline{b}$, we have $|Z|=(k-1)n$. It then follows that $|X|=|Y|=(k-1)n/2$. Thus, if $n$ is odd, we have a contradiction. Assume from now on that $n$ is even. 

Now, consider the sum $S=\sum_{ab\in X}(a+b)+\sum_{\overline{a}\overline{b}\in Y}(a+b)+\sum_{a\overline{b}\in Z}(a+b)$, noting that this is summing over all red edges in $G$. On the one hand, each vertex is included $(2k-2)$ times, so $S=(2k-2)(1+2+\cdots+n+1+2+\cdots+n)=(2k-2)n(n+1)$. On the other hand, for every $i\in[n]$, both $E_i'$ and $F_i'$ contain exactly $k-1$ red edges, each with vertex sum $i$, unless the edge is of the form $ab$, in which case the sum is $i+1$. Thus, $S=2(k-1)(1+2+\cdots n)+|X|=(k-1)n(n+1)+(k-1)n/2$. Therefore, we should have $(2k-2)n(n+1)\equiv(k-1)n(n+1)+(k-1)n/2\pmod n$. However, $(2k-2)n(n+1)-(k-1)n(n+1)-(k-1)n/2=(k-1)n(n+1/2)$ is not a multiple of $n$ since $k-1$ is odd, a contradiction. 
\end{proof}

\subsection{Cyclic nested matchings versus nested matchings}

The goal of this subsection is to confirm Conjecture~4.30 from~\cite{Bavsic2026SomeRO}. Theorem~4.29 from \cite{Bavsic2026SomeRO} already established the exact value of $\rc(\NM_{2k},\NM_{2\ell})$ for all $k,\ell\geq 1$ that are not both even, and Conjecture~4.30 concerns the case when $k,\ell$ are both even. First, we provide a lower bound construction in which, with only a few exceptions, each edge $ij$ is coloured according to $d(i,j)$.

\begin{proposition}\label{prop:nmlower}
For any integers $k,\ell\geq 1$, we have
$$\rc(\NM_{4k},\NM_{4\ell})\geq4k+4\ell-3.$$
\end{proposition}
\begin{proof}
Let $N=4k+4\ell-4$. It suffices to construct a red/blue coloured $N$-vertex cyclically ordered complete graph $G$ that contains no red cyclic copy of $\NM_{4k}$ and no blue cyclic copy of $\NM_{4\ell}$. For every edge $ij$, colour $ij$ red if $1\leq d(i,j)\leq2k-2$, and colour $ij$ blue if $2k-1\leq d(i,j)\leq 2k+2\ell-4$ or if $d(i,j)=2k+2\ell-2$. This assigns a colour to all edges $ij$ in $G$ except those with $d(i,j)=2k+2\ell-3$. For every $i\in[2k+2\ell-2]$, colour the edge connecting $i$ and $i+2k+2\ell-3$ red, and colour all other distance $2k+2\ell-3$ edges blue. 

Suppose for a contradiction that $G$ contains a red cyclic copy of $\NM_{4k}$. Note that there is at most one edge in this copy with distance $2k+2\ell-3$, as any two distance $2k+2\ell-3$ red edges in $G$ intersect. By Lemma~\ref{lemma:nmdistance}\ref{item:1} and~\ref{item:2}, for every $s\in[k-1]$, at most two edges in this copy have their distances in the set $\{2s-1,2s\}$. It follows that the total number of occurrences of the possible red distances $1,2,\ldots,2k-2,2k+2\ell-3$ in this copy is at most $2(k-1)+1=2k-1<2k$, a contradiction. 

Similarly, suppose for a contradiction that $G$ contains a blue cyclic copy of $\NM_{4\ell}$. Note that there do not exist two non-crossing blue edges whose distances lie in $\{2k+2\ell-3, 2k+2\ell-2\}$. Hence, the distances $2k+2\ell-3$ and $2k+2\ell-2$ together can appear at most once in this sequence. By Lemma~\ref{lemma:nmdistance}\ref{item:1} and~\ref{item:2} again, for every $k\leq s\leq k+\ell-2$, at most two edges have their distances in the set $\{2s-1,2s\}$, so the total number of occurrences of the possible blue distances $2k-1,2k,\ldots,2k+2\ell-3,2k+2\ell-2$ in the sequence $d_1,\ldots,d_{2\ell}$ is at most $2(\ell-1)+1=2\ell-1<2\ell$, a contradiction. 
\end{proof}

It remains to prove the matching upper bound.

\begin{theorem}\label{thm:cycnmnm}
    For any positive integers $k,\ell$, we have $$\rc(\NM_{4k},\NM_{4\ell})=4k+4\ell-3.$$
\end{theorem}

\begin{proof}
    Let $N=4k+4\ell-3$. By Proposition~\ref{prop:nmlower}, $\rc(\NM_{4k},\NM_{4\ell})\geq N$, so it suffices to prove the reverse inequality. That is, for any red/blue coloured $N$-vertex cyclically ordered complete graph, there is a red cyclic copy of $\NM_{4k}$ or a blue cyclic copy of $\NM_{4\ell}$.
    Suppose for a contradiction that there is a red/blue edge colouring of the $N$-vertex cyclically ordered complete graph that avoids both a red cyclic $\NM_{4k}$ and a blue cyclic $\NM_{4\ell}$.
    
    Note that both $2k-1$ and $2\ell-1$ are odd numbers.
    Recall from Proposition~\ref{prop:nmcycturan} the notation $E_i=\{\{i-j,i+j\}:1\leq j\leq (N-1)/2\}$ for every $i\in[N]$ and that each $E_i$ is a nested matching. Observe that each set $E_i$ contains exactly one edge of each of the lengths $1,2,\ldots,(N-1)/2$ and no other edges.
    
    Since each $E_i$ contains $2k+2\ell-2$ edges, it must contain exactly $2k-1$ red edges and $2\ell-1$ blue edges, as otherwise it would contain a red copy of $\NM_{4k}$ or a blue copy of $\NM_{4\ell}$.
    Let $c$ be the colour of $\{i-1,i+1\}$, without loss of generality let it be red. Then, we claim that the edges $\{i,i-1\}$ and $\{i,i+1\}$, which lie, respectively, in $E_{i+(N-1)/2}$ and $E_{i-(N-1)/2}$, are also coloured with $c$. Indeed, as $E_i\setminus \{i-1,i+1\}$ contains a blue $\NM_{4\ell-2}$, a blue edge among $\{\{i,i-1\},\{i,i+1\}\}$ would lead to a blue $\NM_{4\ell}$ in $E_i \cup \{\{i,i-1\},\{i,i+1\}\}$.

    We now prove by an induction on $j=1,\ldots,(N-1)/2$ that the two edges $\{i-(j-1),i+j\}$ and $\{i-j,i+(j-1)\}$ have the same colour as $\{i-j,i+j\}$. Indeed, the base case when $j=1$ was established in the previous paragraph. Suppose then that for some $1\leq j\leq (N-3)/2$, for all $1\leq m\leq j$, we have that the two edges $\{i-(m-1),i+m\}$ and $\{i-m,i+(m-1)\}$ have the same colour as $\{i-m,i+m\}$. Assume that $\{i-(j+1),i+(j+1)\}$ has colour $c$, without loss of generality let it be red. We want to show that both $\{i-j,i+j+1\}$ and $\{i-(j+1),i+j\}$ are red as well. 
    
    Suppose for a contradiction, and without loss of generality, that $\{i-j,i+j+1\}$ is blue. From above, $E_i$ contains exactly $2\ell-1$ blue edges. Let $B\subset[(N-1)/2]$ be the set of indices $m$ such that $\{i-m,i+m\}$ is blue, and let $B'=B\cap[j]$. Then, from induction hypothesis, for every $m\in B'$, $\{i-(m-1),i+m\}$ is blue as well. This implies that 
    \[\left(\bigcup_{m\in B'}\{i-(m-1),i+m\}\right)\cup\{i-j,i+j+1\}\cup\left(\bigcup_{m\in B\setminus B'}\{i-m,i+m\}\right)\]
    forms a blue cyclic copy of $\NM_{4\ell}$, a contradiction. Therefore, both of those edges are red. See Figure~\ref{fig:nmvsnm} for an illustration.

    \begin{figure}[H]
\begin{subfigure}[c]{0.48\textwidth}
\centering
\scalebox{0.85}{%
\begin{tikzpicture}
\tikzstyle{every node}=[font=\tiny]
\node[draw=none,minimum size=6cm,regular polygon,regular polygon sides=19] (a) {};

\node[draw=none,minimum size=7.5cm,regular polygon,regular polygon sides=19] (b) {};
\node[draw=none,minimum size=6.75cm,regular polygon,regular polygon sides=19] (c) {};

\foreach \x in {1,2,3,5,6,7,9,10,11,12,14,15,16,18,19}
  \fill (a.corner \x) circle[radius=1pt];

\draw (c.corner 1) node {$i$};
\draw (c.corner 2) node {$i-1$};
\draw (c.corner 19) node {$i+1$};
\draw (c.corner 3) node {$i-2$};
\draw (c.corner 18) node {$i+2$};
\draw (b.corner 5) node {$i-(j-1)$};
\draw (b.corner 16) node {$i+(j-1)$};
\draw (b.corner 6) node {$i-j$};
\draw (b.corner 15) node {$i+j$};
\draw (b.corner 7) node {$i-(j+1)$};
\draw (b.corner 14) node {$i+(j+1)$};
\draw (b.corner 9) node {$i-\frac{N-3}{2}$};
\draw (b.corner 12) node {$i+\frac{N-3}{2}$};
\draw (b.corner 10) node {$i-\frac{N-1}{2}$};
\draw (b.corner 11) node {$i+\frac{N-1}{2}$};

\draw (a.corner 4) node {$\cdots$};
\draw (a.corner 8) node {$\cdots$};
\draw (a.corner 13) node {$\cdots$};
\draw (a.corner 17) node {$\cdots$};

\draw[red] (a.corner 2)--(a.corner 19);
\draw[blue] (a.corner 3)--(a.corner 18);
\draw[red] (a.corner 5)--(a.corner 16);
\draw[blue] (a.corner 6)--(a.corner 15);
\draw[red] (a.corner 7)--(a.corner 14);
\draw[blue] (a.corner 9)--(a.corner 12);
\draw[red] (a.corner 10)--(a.corner 11);

\draw[red] (a.corner 1)--(a.corner 2);
\draw[red] (a.corner 1)--(a.corner 19);
\draw[blue] (a.corner 2)--(a.corner 18);
\draw[blue] (a.corner 19)--(a.corner 3);
\draw[blue] (a.corner 5)--(a.corner 15);
\draw[blue] (a.corner 6)--(a.corner 16);
\draw[red,dashed] (a.corner 7)--(a.corner 15);
\draw[red,dashed] (a.corner 6)--(a.corner 14);
\end{tikzpicture}
}
\end{subfigure}
\begin{subfigure}[c]{0.48\textwidth}
\centering
\scalebox{0.85}{%
\begin{tikzpicture}
\tikzstyle{every node}=[font=\tiny]
\node[draw=none,minimum size=6cm,regular polygon,regular polygon sides=19] (a) {};

\node[draw=none,minimum size=7.5cm,regular polygon,regular polygon sides=19] (b) {};
\node[draw=none,minimum size=6.75cm,regular polygon,regular polygon sides=19] (c) {};

\foreach \x in {1,2,3,5,6,7,9,10,11,12,14,15,16,18,19}
  \fill (a.corner \x) circle[radius=1pt];

\draw (c.corner 1) node {$i$};
\draw (c.corner 2) node {$i-1$};
\draw (c.corner 19) node {$i+1$};
\draw (c.corner 3) node {$i-2$};
\draw (c.corner 18) node {$i+2$};
\draw (b.corner 5) node {$i-(j-1)$};
\draw (b.corner 16) node {$i+(j-1)$};
\draw (b.corner 6) node {$i-j$};
\draw (b.corner 15) node {$i+j$};
\draw (b.corner 7) node {$i-(j+1)$};
\draw (b.corner 14) node {$i+(j+1)$};
\draw (b.corner 9) node {$i-\frac{N-3}{2}$};
\draw (b.corner 12) node {$i+\frac{N-3}{2}$};
\draw (b.corner 10) node {$i-\frac{N-1}{2}$};
\draw (b.corner 11) node {$i+\frac{N-1}{2}$};

\draw (a.corner 4) node {$\cdots$};
\draw (a.corner 8) node {$\cdots$};
\draw (a.corner 13) node {$\cdots$};
\draw (a.corner 17) node {$\cdots$};

\draw[ultra thick, blue] (a.corner 2)--(a.corner 18);
\draw[ultra thick, blue] (a.corner 5)--(a.corner 15);
\draw[ultra thick, blue] (a.corner 6)--(a.corner 14);
\draw[ultra thick, blue] (a.corner 9)--(a.corner 12);

\draw[red] (a.corner 2)--(a.corner 19);
\draw[blue] (a.corner 3)--(a.corner 18);
\draw[red] (a.corner 5)--(a.corner 16);
\draw[blue] (a.corner 6)--(a.corner 15);
\draw[red] (a.corner 7)--(a.corner 14);
\draw[blue] (a.corner 9)--(a.corner 12);
\draw[red] (a.corner 10)--(a.corner 11);

\draw[red] (a.corner 1)--(a.corner 2);
\draw[red] (a.corner 1)--(a.corner 19);
\draw[blue] (a.corner 2)--(a.corner 18);
\draw[blue] (a.corner 19)--(a.corner 3);
\draw[blue] (a.corner 5)--(a.corner 15);
\draw[blue] (a.corner 6)--(a.corner 16);
\draw[red,dashed] (a.corner 7)--(a.corner 15);
\draw[blue] (a.corner 6)--(a.corner 14);

\end{tikzpicture}
}
\end{subfigure}

\caption{An illustration of the proof of Theorem~\ref{thm:cycnmnm}. Assuming for a contradiction that there is a red/blue colouring of $K_N^{\mathrm{cyc}}$ for $N=4k+4\ell -3$ with no red $\NM_{4k}$ and no blue $\NM_{4 \ell}$. The set $E_i$ contains $2k-1$ red edges and $2\ell -1$ blue edges. In the induction step, assuming that $\{i-(j+1),i+(j+1)\}$ is red, then the edges $\{i-j,i+(j+1)\}$ and $\{i-(j+1),i+j\}$ (marked by dashed lines) also have to be red. Indeed, if one of them is blue, then we can find a blue cyclic copy of $\NM_{4\ell}$, as depicted on the right.}\label{fig:nmvsnm}
\end{figure}
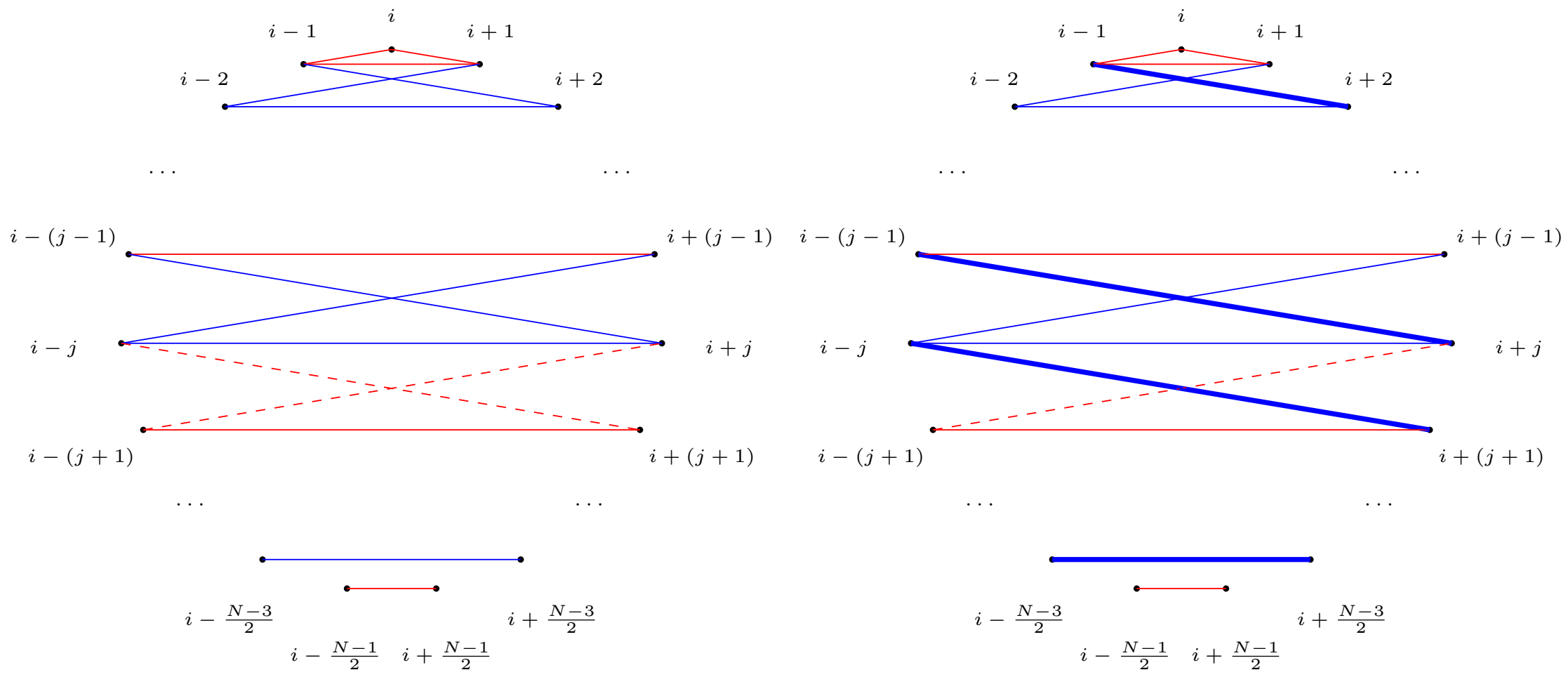

    At the end of the induction, we get that the colouring of $E_i$ uniquely determines the colouring of $E_{i+(N-1)/2}$ and $E_{i-(N-1)/2}$. Repeating this argument, we have that the colouring of $E_{i+(N-1)/2}$ determines the colouring of $E_{i-1}$, so the colouring of $E_i$ uniquely determines the colouring of $E_{i-1}$. More specifically, for every $j\in[(N-1)/2]$, the colour of $\{i-j,i+j\}$ in $E_i$ is the same as the colour of $\{i-j,i+j-1\}$ in $E_{i+(N-1)/2}$, which is the same as the colour of $\{i-j-1,i+j-1\}$ in $E_{i-1}$. Since this is true for all $i\in[N]$ and $j\in[(N-1)/2]$, we get that in particular, all edges in $G$ with the same distance have the same colour. Furthermore, for any fixed $i$, the induction above implies that for every $j\in[(N-1)/4]$, the edges $\{i-j,i+j\}$ and $\{i-j,i+j-1\}$ have the same colour. Thus, the common colours in $G$ of all distance $2j$ and distance $2j-1$ edges are the same.

    However, recall that $E_i$ consists of $(N-1)/2$ edges, with each distance in $[(N-1)/2]$ appearing exactly once.
    Therefore, the edges of $E_i$ can be partitioned into $(N-1)/4$ pairs of edges of length $2j-1$ and $2j$ for $j\in[(N-1)/4]$, with edges in the same pair having the same colour.
    This implies that $E_i$ contains an even number of red edges and an even number of blue edges, contradicting that $E_i$ contains $2k-1$ red edges and $2 \ell-1$ blue edges. 
\end{proof}

Combining this result with Theorem 4.29 from \cite{Bavsic2026SomeRO} we obtain the following formula for $\rc(\NM_{2k},\NM_{2\ell})$.

\begin{corollary}
    Let $k,\ell \geq 1$. Then
    $$\rc(\NM_{2k},\NM_{2\ell})=2k+2\ell-2 - \mathbbm{1}\{(2 \mid k) \land (2 \mid \ell)\}.$$
\end{corollary}

\subsection{Ordered nested matchings versus monotone paths}

In this subsection, we confirm Conjecture~4.31 from~\cite{Bavsic2026SomeRO}. For the upper bound, we use the following result which explicitly appeared in Balko, Cibulka, Král, and Kynčl \cite{balko2013ramsey}, and which was already implicit in the proof of Theorem~4.3 in Károlyi, Pach, and Tóth \cite{karolyi1997ramsey}.

\begin{lemma}[{\cite[Lemma 18]{balko2013ramsey}}]\label{lemma:ordkkmp}
For all positive integers $k,\ell$, we have
$$\ro(K_k,\MP_\ell)=1+(k-1)(\ell-1).$$
\end{lemma}

We now confirm Conjecture~4.31 from \cite{Bavsic2026SomeRO} by using a lower bound construction in which, again, each edge $ij$ is coloured according to $d(i,j)$.
\begin{theorem}\label{thm:ordnmmp}
For all positive integers $k,\ell$, we have
$$\ro(\NM_{2k},\MP_\ell)=1+(2k-1)(\ell-1).$$
\end{theorem}
\begin{proof}
Since $\NM_{2k}$ is a subgraph of $K_{2k}$, the upper bound follows from the upper bound part of Lemma~\ref{lemma:ordkkmp}. 

For the lower bound, consider the following red/blue coloured ordered complete graph $G$ on $(2k-1)(\ell-1)$ vertices. For any $1\leq i<j\leq(2k-1)(\ell-1)$, the edge $ij$ is coloured red if $j-i\leq 2k-2$, and blue otherwise. Note that the red graph does not contain an ordered copy of $\NM_{2k}$, as no red edge $ij$ satisfies $j-i\geq 2k-1$. In contrast, any blue edge $ij$ in $G$ satisfies $j-i\geq 2k-1$. Suppose there exists $1\leq v_1<v_2<\cdots<v_\ell\leq(2k-1)(\ell-1)$ such that they form a blue ordered copy of $\MP_\ell$. Then, $v_{i+1}-v_i\geq 2k-1$ for every $i\in[\ell-1]$, so $v_{\ell}-v_1\geq(2k-1)(\ell-1)$, a contradiction. Hence, there is also no blue ordered copy of $\MP_\ell$, proving the lower bound. 
\end{proof}

\subsection{Cyclic nested matchings versus monotone paths}

In this subsection, we prove the lower bound part of Conjecture~4.32 as follows. Once again, each edge $ij$ is coloured almost entirely according to $d(i,j)$.

\begin{proposition}\label{prop:lower}
For any $k\geq 1$ and $\ell\geq3$, 
\[\rc(\NM_{2k},\MP_\ell)\geq\begin{cases}k\ell-k+1,&\text{if }k\text{ is odd,}\\k\ell-k-\ell+3,&\text{if }k\text{ is even.}\end{cases}\]
\end{proposition}
\begin{proof}
First suppose that $k$ is odd. Let $N=k\ell-k$, and consider the red/blue coloured cyclically ordered $N$-vertex complete graph $G$ in which an edge $ij$ is coloured red if $d(i,j)\leq k-1$, and blue otherwise. Then, there is no red cyclic copy of $\NM_{2k}$ in $G$ as any such copy must contain an edge $ij$ with $d(i,j)\geq k$ by Lemma~\ref{lemma:nmdistance}. Suppose for a contradiction that there is a blue cyclic copy of $\MP_\ell$ in $G$ with vertices $v_1,\ldots,v_\ell$ in clockwise order, such that $v_iv_{i+1}$ is an edge for every $i\in[\ell-1]$. By symmetry, we may assume that $v_1=1$. Then, $d_c(v_i,v_{i+1})\geq k$ for every $i\in[\ell-1]$, so $d_c(v_1,v_\ell)\geq k(\ell-1)=N$, a contradiction. Thus, $\rc(\NM_{2k},\MP_\ell)\geq k\ell-k+1$.

Now suppose that $k$ is even, and let $N=k\ell-k-\ell+2$. Consider the following red/blue coloured ordered $N$-vertex complete graph $G$. For every $1\leq i<j\leq N$, colour edge $ij$ red if $d(i,j)\leq k-2$, or if $j-i=k-1$ and $i\in[k]$. All other edges are coloured blue. 

Suppose for a contradiction that $G$ contains a red cyclic copy of $\NM_{2k}$. The distance $k-1$ can appear at most once in this sequence, as no two red edges with distance $k-1$ are non-crossing. Furthermore, by Lemma~\ref{lemma:nmdistance}, for every $1\leq s\leq(k-2)/2$, at most two edges in this copy have distances in $\{2s-1,2s\}$. It follows that the total number of occurrences of the allowed red distances $1,2,\ldots,k-1$ in the sequence $d_1,\ldots,d_{k}$ is at most $2(k-2)/2+1=k-1<k$, a contradiction. 

Now suppose for a contradiction that $G$ contains a blue cyclic copy of $\MP_\ell$ with vertices $v_1,\ldots,v_\ell$ in clockwise order, such that $v_iv_{i+1}$ is an edge for every $i\in[\ell-1]$. Then, $d_c(v_1,v_\ell)\leq N-1=k\ell-k-\ell+1=(k-1)(\ell-1)$. Since $d_c(v_1,v_\ell)=\sum_{i=1}^{\ell-1}d_c(v_i,v_{i+1})$ and $d_c(v_i,v_{i+1})\geq k-1$ for every $i\in[\ell-1]$, we must have $d_c(v_i,v_{i+1})=k-1$ for every $i\in[\ell-1]$. However, $\{v_1,\ldots,v_{\ell-1}\}\cap[k]=\varnothing$, as by definition the edge between $i$ and $i+k-1$ is red for every $i\in[k]$. Therefore, $(\ell-2)(k-1)=d_c(v_1,v_{\ell-1})\leq N-k-1=(\ell-2)(k-1)-1$, a contradiction. 
\end{proof}

\section{Cyclic nested matchings versus complete graphs}\label{sec:disprove}
Surprisingly, Conjecture~4.36 from \cite{Bavsic2026SomeRO} turns out to be false. In the case when $\ell=3$, a construction of Cs\'oka, Bl\'azsik, Kir\'aly, and Lenger~{\cite[Theorem 3, Construction 2]{MR4442957}} in the context of necklace folding problems already implies that $\rc(\NM_{2k},K_3)\geq(2+\sqrt 2+o(1))k>2k$ for infinitely many $k$. Generalizing their construction, we obtain counterexamples to~\cite[Conjecture 4.36]{Bavsic2026SomeRO} for every fixed $\ell\geq 3$ and arbitrarily large $k$.  

\begin{theorem}\label{thm:disproof}
For every fixed $\ell\geq3$, there exists arbitrarily large $k$, such that \[\rc(\NM_{2k},K_\ell)\geq(\ell-1+\sqrt{\ell-1}+o(1))k>k\ell.\]
\end{theorem}

We will prove Theorem~\ref{thm:disproof} via the following result.
\begin{proposition}\label{prop:disproof}
For every fixed $r\geq2$, there exist arbitrarily large $N$ and an $N$-vertex cyclically ordered graph $H$, such that the following hold. The complement of $H$ is $r$-partite, and $H$ does not contain a cyclic nested matching with $\left(\frac1{r+\sqrt r}+o(1)\right)N$ edges.
\end{proposition}

Assuming Proposition~\ref{prop:disproof}, Theorem~\ref{thm:disproof} follows immediately.
\begin{proof}[Proof of Theorem~\ref{thm:disproof}]
Let $\ell\geq3$ be fixed. By Proposition~\ref{prop:disproof}, there exist arbitrarily large $N$, such that there is an $N$-vertex ordered graph $H$ that does not contain a cyclic nested matching with $\left(\frac1{\ell-1+\sqrt{\ell-1}}+o(1)\right)N$ edges, and whose complement is $(\ell-1)$-partite. Let $G$ be the red/blue coloured $N$-vertex cyclically ordered complete graph, such that the red subgraph of $G$ is $H$, and the blue subgraph of $G$ is the complement of $H$. Then, there is no blue cyclic copy of $K_\ell$, as the blue graph is $(\ell-1)$-partite. Let $k$ be the number of edges in the largest red cyclic nested matching in $G$. Then, $k\leq\left(\frac1{\ell-1+\sqrt{\ell-1}}+o(1)\right)N$ from assumption, so $N\geq(\ell-1+\sqrt{\ell-1}+o(1))k$. Thus, $\rc(\NM_{2k},K_\ell)\geq\left(\ell-1+\sqrt{\ell-1}+o(1)\right)k>k\ell$.  
\end{proof}

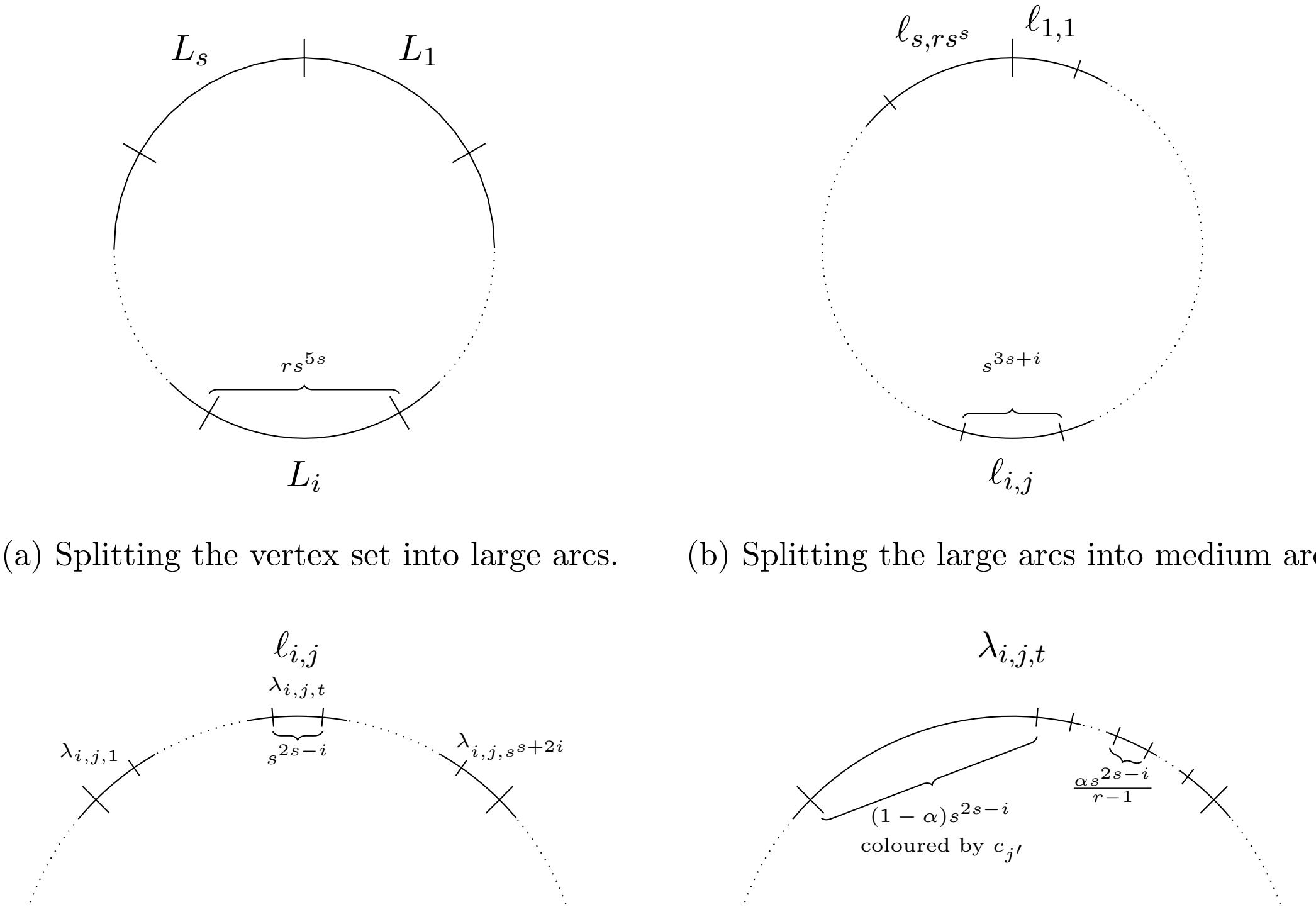
\begin{figure}[h]
\centering
\begin{subfigure}[c]{0.48\textwidth}
\centering
\scalebox{1}{%
\begin{tikzpicture}
\draw [domain=0:180] plot ({2*cos(\x)}, {2*sin(\x)});
\draw [dotted,domain=180:225] plot ({2*cos(\x)}, {2*sin(\x)});
\draw [domain=225:315] plot ({2*cos(\x)}, {2*sin(\x)});
\draw [dotted,domain=315:360] plot ({2*cos(\x)}, {2*sin(\x)});
\draw (0,1.8)--(0,2.2);
\draw [rotate=60] (0,1.8)--(0,2.2);
\draw [rotate=-60] (0,1.8)--(0,2.2);
\draw [rotate=210] (0,1.8)--(0,2.2);
\draw [rotate=150] (0,1.8)--(0,2.2);
\draw [rotate=30] node at (0,2.4) {$L_s$};
\draw [rotate=-30] node at (0,2.4) {$L_1$};
\draw [rotate=180] node at (0,2.4) {$L_i$};
\draw [decorate, decoration={brace, mirror}] ({2*cos(300)}, {2*sin(300)+0.2})--({2*cos(240)}, {2*sin(240)+0.2});
\draw node at (0,-1.2) {\tiny{$rs^{5s}$}};

\draw node at (0,2.7) {};
\draw node at (0,-2.7) {};
\end{tikzpicture}
}
\caption{Splitting the vertex set into large arcs.}
\end{subfigure}
\begin{subfigure}[c]{0.48\textwidth}
\centering
\scalebox{1.0}{%
\begin{tikzpicture}
\draw [domain=60:140] plot ({2*cos(\x)}, {2*sin(\x)});
\draw [dotted,domain=140:245] plot ({2*cos(\x)}, {2*sin(\x)});
\draw [domain=245:295] plot ({2*cos(\x)}, {2*sin(\x)});
\draw [dotted,domain=295:420] plot ({2*cos(\x)}, {2*sin(\x)});
\draw (0,1.8)--(0,2.2);
\draw [rotate=40] (0,1.9)--(0,2.1);
\draw [rotate=-20] (0,1.9)--(0,2.1);
\draw [rotate=195] (0,1.9)--(0,2.1);
\draw [rotate=165] (0,1.9)--(0,2.1);
\draw [rotate=20] node at (0,2.4) {$\ell_{s,rs^s}$};
\draw [rotate=-10] node at (0,2.4) {$\ell_{1,1}$};
\draw [rotate=180] node at (0,2.4) {$\ell_{i,j}$};
\draw [decorate, decoration={brace, mirror}] ({2*cos(285)}, {2*sin(285)+0.15})--({2*cos(255)}, {2*sin(255)+0.15});
\draw node at (0,-1.2) {\tiny{$s^{3s+i}$}};

\draw node at (0,2.7) {};
\draw node at (0,-2.7) {};
\end{tikzpicture}
}
\caption{Splitting the large arcs into medium arcs.}
\end{subfigure}
\begin{subfigure}[c]{0.48\textwidth}
\centering
\scalebox{1}{%
\begin{tikzpicture}
\draw [domain=40:60] plot ({3*cos(\x)}, {3*sin(\x)});
\draw [domain=120:140] plot ({3*cos(\x)}, {3*sin(\x)});
\draw [domain=80:100] plot ({3*cos(\x)}, {3*sin(\x)});
\draw [dotted,domain=140:160] plot ({3*cos(\x)}, {3*sin(\x)});
\draw [dotted,domain=20:40] plot ({3*cos(\x)}, {3*sin(\x)});
\draw [dotted,domain=100:120] plot ({3*cos(\x)}, {3*sin(\x)});
\draw [dotted,domain=60:80] plot ({3*cos(\x)}, {3*sin(\x)});

\draw [rotate=45] (0,2.8)--(0,3.2);
\draw [rotate=-45] (0,2.8)--(0,3.2);
\draw [rotate=5] (0,2.9)--(0,3.1);
\draw [rotate=-5] (0,2.9)--(0,3.1);
\draw [rotate=35] (0,2.9)--(0,3.1);
\draw [rotate=-35] (0,2.9)--(0,3.1);

\draw node at (0,3.3) {\tiny{$\lambda_{i,j,t}$}};
\draw [rotate=40] node at (0,3.4) {\tiny{$\lambda_{i,j,1}$}};
\draw [rotate=-40] node at (0,3.5) {\tiny{$\lambda_{i,j,s^{s+2i}}$}};
\draw node at (0,3.7) {$\ell_{i,j}$};

\draw [decorate, decoration={brace}] ({3*cos(85)}, {3*sin(85)-0.15})--({3*cos(95)}, {3*sin(95)-0.15});
\draw node at (0,2.6) {\tiny{$s^{2s-i}$}};

\draw node at (0,4.3) {};
\end{tikzpicture}
}
\caption{Splitting medium arcs into small arcs.}
\end{subfigure}
\begin{subfigure}[c]{0.48\textwidth}
\centering
\scalebox{1.0}{%
\begin{tikzpicture}
\draw [domain=40:54] plot ({3*cos(\x)}, {3*sin(\x)});
\draw [domain=76:140] plot ({3*cos(\x)}, {3*sin(\x)});
\draw [domain=59.5:70.5] plot ({3*cos(\x)}, {3*sin(\x)});
\draw [dotted,domain=140:160] plot ({3*cos(\x)}, {3*sin(\x)});
\draw [dotted,domain=20:40] plot ({3*cos(\x)}, {3*sin(\x)});
\draw [dotted,domain=70.5:76] plot ({3*cos(\x)}, {3*sin(\x)});
\draw [dotted,domain=54:59.5] plot ({3*cos(\x)}, {3*sin(\x)});

\draw [rotate=45] (0,2.8)--(0,3.2);
\draw [rotate=-45] (0,2.8)--(0,3.2);
\draw [rotate=-5] (0,2.9)--(0,3.1);
\draw [rotate=-12] (0,2.9)--(0,3.1);
\draw [rotate=-38] (0,2.9)--(0,3.1);
\draw [rotate=-21.5] (0,2.9)--(0,3.1);
\draw [rotate=-28.5] (0,2.9)--(0,3.1);

\draw node at (0,3.7) {$\lambda_{i,j,t}$};

\draw [decorate, decoration={brace}] ({3*cos(63)}, {3*sin(63)-0.15})--({3*cos(70)}, {3*sin(70)-0.15});
\draw [rotate=-25] node at (0,2.5) {\tiny{$\frac{\alpha s^{2s-i}}{r-1}$}};

\draw [decorate, decoration={brace}] ({3*cos(85)}, {3*sin(85)-0.2})--({3*cos(132)}, {3*sin(132)-0.3});
\draw [rotate=21] node at (0,2.1) {\tiny{$(1-\alpha) s^{2s-i}$}};
\draw [rotate=24.5] node at (0,1.75) {\tiny{coloured by $c_{j'}$}};

\draw node at (0,4.3) {};
\end{tikzpicture}
}
\caption{Splitting small arcs.}
\end{subfigure}
\caption{An illustration of the colouring of the vertices in Construction~\ref{construction}.}
\label{fig:construction}
\end{figure}

The proof of Proposition~\ref{prop:disproof} is quite technical. We first describe the construction of $H$.
\begin{construction}\label{construction}
Fix $r\geq 2$, and let $\alpha=1-1/\sqrt r$. For every $s\geq 2$, let $n=s^{5s+1}$, and let $N=rs^{5s+1}$. Define an $N$-vertex ordered graph $H$ with vertex set $[N]$ as follows. Each vertex of $H$ will be coloured with a colour $c_1,\ldots,c_r$, so that there are exactly $s^{5s+1}$ vertices with colour $c_i$ for every $i\in[r]$. 

To begin, for every $i\in[s]$, let $L_i$ be the $rs^{5s}$ vertices $(i-1)rs^{5s}+1,(i-1)rs^{5s}+2,\ldots,irs^{5s}$. Call $L_1,L_2,\ldots,L_s$ the \emph{large arcs}.

Then, for each large arc $L_i$, subdivide $L_i$ into $rs^{2s-i}$ \emph{medium arcs} $\ell_{i,1},\ell_{i,2},\ldots,\ell_{i,rs^{2s-i}}$, each containing $s^{3s+i}$ consecutive vertices. For every $j\in[rs^{2s-i}]$, let $j'\in[r]$ be the unique index satisfying $j\equiv j'\pmod r$, and call $c_{j'}$ the \emph{majority colour} of the medium arc $\ell_{i,j}$.

Next, we colour the vertices on each medium arc as follows. For every $i\in[s]$ and every $j\in[rs^{2s-i}]$, suppose $c_{j'}$ is the majority colour of the medium arc $\ell_{i,j}$. Subdivide $\ell_{i,j}$ into $s^{s+2i}$ \emph{small arcs} $\lambda_{i,j,1},\lambda_{i,j,2},\ldots,\lambda_{i,j,s^{s+2i}}$, each containing $s^{2s-i}$ consecutive vertices. For every $t\in[s^{s+2i}]$, colour the first $(1-\alpha)s^{2s-i}$ vertices on the small arc $\lambda_{i,j,t}$ with the majority colour $c_{j'}$, and partition the last $\alpha s^{2s-i}$ vertices on $\lambda_{i,j,t}$ into $r-1$ \emph{tiny arcs} of $\alpha s^{2s-i}/(r-1)$ consecutive vertices. Then, arbitrarily assign to each tiny arc a distinct colour that is not $c_{j'}$, and colour all vertices in each tiny arc with the assigned colour. Note that in this way, the majority colour $c_{j'}$ is the colour that appears the most often within each small arc $\lambda_{i,j,t}$. See Figure~\ref{fig:construction} for an illustration.

Finally, we define $H$ to be the ordered graph with vertex set $[N]$, in which two vertices $i$ and $j$ form an edge in $H$ if and only if they have the same colour. Then, the complement of $H$ is evidently $r$-partite. 
\end{construction}

\begin{proof}[Proof of Proposition~\ref{prop:disproof}]
Fix $r\geq 2$, and let $\alpha=1-1/\sqrt r$. For every $s\geq 2$, let $n=s^{5s+1}$, and let $N=rs^{5s+1}$. Let $H$ be the $N$-vertex ordered graph given by Construction~\ref{construction}. As usual, we view the vertices of $H$ as being arranged on a circle. Let $M$ be a maximum cyclic nested matching in $H$. 

From the definition of a cyclic nested matching, there is a chord $C$ of the circle that splits the circle into two halves, such that every edge in $M$ contains exactly one vertex in each half. By rotating $C$ slightly if necessary, we may assume that each end of this chord $C$ lies on a large arc. Let $1\leq a\leq b\leq s$ be the indices such that one end of $C$ lies on $L_a$, and the other end lies on $L_b$. If $a=b$, then one end of each edge in $M$ is in $L_a$, so $|M|\leq|L_a|=N/s=o(N)$, as required. Thus, we may assume $a<b$. Let $M'\subset M$ be the set of edges in $M$ containing no vertices in $L_a\cup L_b$. Then, $|M'|\geq|M|-|L_a|-|L_b|=|M|-2rs^{5s}=|M|-2N/s=|M|-o(N)$ as $s$ tends to infinity.

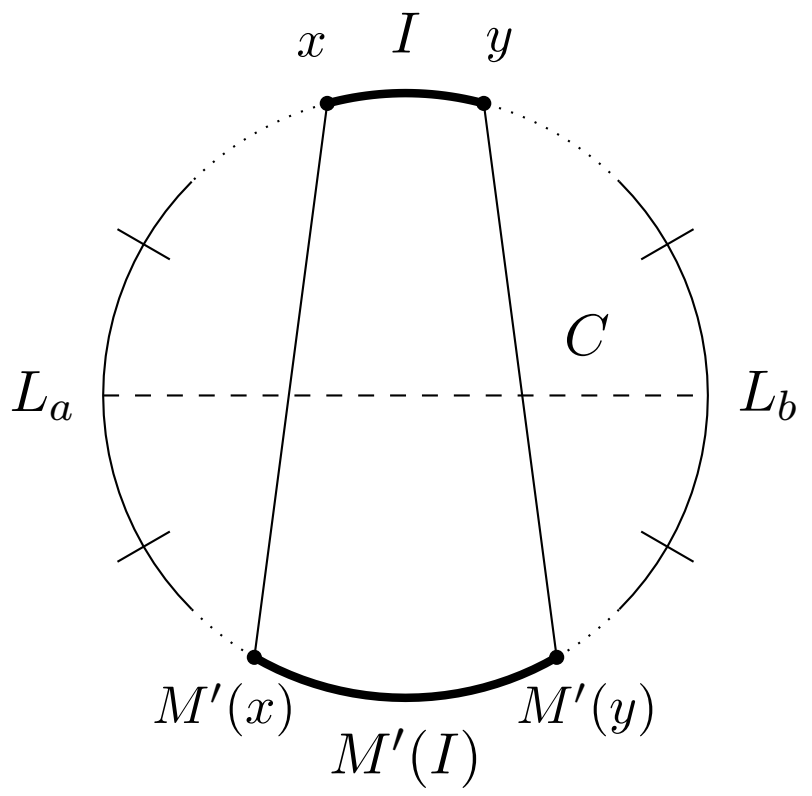
\begin{figure}[h]
\centering
\scalebox{1}{%
\begin{tikzpicture}
\draw [domain=-45:45] plot ({2*cos(\x)}, {2*sin(\x)});
\draw [domain=135:225] plot ({2*cos(\x)}, {2*sin(\x)});
\draw [ultra thick, domain=75:105] plot ({2*cos(\x)}, {2*sin(\x)});
\draw [ultra thick, domain=-60:-120] plot ({2*cos(\x)}, {2*sin(\x)});
\draw [dotted, domain=45:75] plot ({2*cos(\x)}, {2*sin(\x)});
\draw [dotted, domain=105:135] plot ({2*cos(\x)}, {2*sin(\x)});
\draw [dotted, domain=225:240] plot ({2*cos(\x)}, {2*sin(\x)});
\draw [dotted, domain=300:315] plot ({2*cos(\x)}, {2*sin(\x)});
\fill ({2*cos(105)}, {2*sin(105}) circle [radius=0.05];
\fill ({2*cos(75)}, {2*sin(75}) circle [radius=0.05];
\fill ({2*cos(-60)}, {2*sin(-60}) circle [radius=0.05];
\fill ({2*cos(-120)}, {2*sin(-120}) circle [radius=0.05];
\draw [rotate=-60] (0,1.8)--(0,2.2);
\draw [rotate=-120] (0,1.8)--(0,2.2);
\draw [rotate=120] (0,1.8)--(0,2.2);
\draw [rotate=60] (0,1.8)--(0,2.2);
\draw [dashed] (-2,0)--(2,0);
\draw ({2*cos(105)}, {2*sin(105})--({2*cos(-120)}, {2*sin(-120});
\draw ({2*cos(75)}, {2*sin(75})--({2*cos(-60)}, {2*sin(-60});
\draw node at (1.2,0.4) {$C$};
\draw node at (0,2.4) {$I$};
\draw [rotate=180] node at (0,2.4) {$M'(I)$};
\draw [rotate=90] node at (0,2.4) {$L_a$};
\draw [rotate=-90] node at (0,2.4) {$L_b$};
\draw [rotate=15] node at (0,2.4) {\small{$x$}};
\draw [rotate=-15] node at (0,2.4) {\small{$y$}};
\draw [rotate=150] node at (0,2.4) {\small{$M'(x)$}};
\draw [rotate=-150] node at (0,2.4) {\small{$M'(y)$}};

\draw node at (0,2.7) {};
\draw node at (0,-2.7) {};
\end{tikzpicture}
}

\caption{An illustration of the definition of $C$ as well as of $M'(I)$ for an interval of consecutive vertices disjoint from $L_a$ and $L_b$. Note that $M'(I)$ is also an interval by its definition.}

\label{fig:mprime}
\end{figure}

For every vertex $x$ on $H$ covered by $M'$, let $M'(x)$ be its unique neighbour in $M'$. For every set $I$ of consecutive vertices on $H$ disjoint from $L_a\cup L_b$, if $x$ and $y$ are the clockwise first and last vertex in $I$ that is covered by $M'$, respectively, then let $M'(I)$ be the set of clockwise consecutive vertices on $H$ from $M'(y)$ to $M'(x)$. See Figure~\ref{fig:mprime} for an illustration.

For any $1\leq i<i'\leq s$, we call $(i,i')$ a \emph{bonded pair} if there exists an edge in $M'$ with one vertex in $L_i$ and the other in $L_{i'}$. We claim that the number of bonded pairs is at most $s-3$. Indeed, form the auxiliary graph $H'$ whose vertex set is $[s]\setminus\{a,b\}$, and whose edge set consists of the bonded pairs. Then, $H'$ must be bipartite from the definition of $a$ and $b$. Moreover, $H'$ must also be acyclic, as the set of edges corresponding to an even cycle in $H'$ would contain two edges that cross each other. Therefore, $H'$ is a forest on $s-2$ vertices, so has at most $s-3$ edges, and thus there are at most $s-3$ bonded pairs. 

Now, let $(i,i')$ be any bonded pair with $i<i'$. Recall from Construction~\ref{construction} that for every $j\in[s^{2s-i}]$, the majority colours of the medium arcs $\ell_{i,(j-1)r+1},\ell_{i,(j-1)r+2},\ldots,\ell_{i,jr}$ are $c_1,c_2,\ldots,c_r$, respectively. Let $\overline{\ell}_{i,j}=\bigcup_{p=1}^{r}\ell_{i,(j-1)r+p}$. We say that an index $j\in[s^{2s-i}]$ is \emph{$(i,i')$-regular}, if there exists $j'\in[rs^{2s-i'}]$ such that $M'(\overline{\ell}_{i,j})\subset\ell_{i',j'}$. An edge in $M'$ is called \emph{regular}, if for some bonded pair $(i,i')$, one of its endpoints is in $\overline{\ell}_{i,j}$ for some $j\in[s^{2s-i}]$ that is $(i,i')$-regular.

Let $M''\subset M'$ be the set of regular edges. We claim that $|M''|\geq|M'|-o(N)$. Indeed, let $(i,i')$ be any bonded pair, and consider the number of indices $j\in[s^{2s-i}]$ that are not $(i,i')$-regular, but satisfy that $M'(\overline{\ell}_{i,j})\cap L_{i'}\not=\varnothing$. Temporarily, let $\ell_{i',0}$ be the medium arc right before $\ell_{i',1}$ in $L_{i'-1}$, and let $\ell_{i',rs^{2s-i'}+1}$ be the medium arc right after $\ell_{i',rs^{2s-i'}}$ in $L_{i'+1}$. Then, from definition, $j$ being not $(i,i')$-regular implies that there is some $0\leq j'\leq rs^{2s-i'}$, such that $M'(\overline{\ell}_{i,j})$ intersects both $\ell_{i',j'}$ and $\ell_{i',j'+1}$. On the other hand, every such $0\leq j'\leq rs^{2s-i'}$ can correspond to at most one $j\in[s^{2s-i}]$ that is not $(i,i')$-regular, so there are at most $rs^{2s-i'}+1$ indices $j\in[s^{2s-i}]$ with $M'(\overline{\ell}_{i,j})\cap L_{i'}\not=\varnothing$ that are not $(i,i')$-regular. Recall that each medium arc contains $s^{3s+i}$ vertices, so at most $(rs^{2s-i'}+1)rs^{3s+i}<2r^2s^{5s+i-i'}\leq2r^2s^{5s-1}$ edges in the bonded pair $(i,i')$ are not regular. Therefore, across the at most $s-3$ bonded pairs, the number of edges in $M'$ that are not regular is at most $(s-3)2r^2s^{5s-1}<2r^2s^{5s}=o(N)$, proving the claim.

From now on we focus on $M''$. Fix any bonded pair $(i,i')$ and any $j\in[s^{2s-i}]$ that is $(i,i')$-regular, and let $j'\in[rs^{2s-i'}]$ satisfy $M'(\overline{\ell}_{i,j})\subset\ell_{i',j'}$.

Define the \emph{efficiency} of $M''$ on $\overline{\ell}_{i,j}$ to be 
\[\eff(\overline{\ell}_{i,j})=\frac{|\{x\in\overline{\ell}_{i,j}\cup M'(\overline{\ell}_{i,j}):x\in V(M'')\}|}{|\overline{\ell}_{i,j}|+|M'(\overline{\ell}_{i,j})|}.\]
In other words, $\eff(\overline{\ell}_{i,j})$ is the proportion of vertices in $\overline{\ell}_{i,j}\cup M'(\overline{\ell}_{i,j})$ that is covered by the matching $M''$. We claim that $\eff(\overline{\ell}_{i,j})\leq\frac2{r+\sqrt r}+o(1)$.

Recall that $\overline{\ell}_{i,j}=\bigcup_{p=1}^{r}\ell_{i,(j-1)r+p}$, and the medium arc $\ell_{i,(j-1)r+p}$ has majority colour $c_p$ for every $p\in[r]$. Recall also that from Construction~\ref{construction}, for every $p\in[r]$, $\ell_{i,(j-1)r+p}$ is further subdivided into $s^{s+2i}$ small arcs $\lambda_{i,(j-1)r+p,t}$ of length $s^{2s-i}$, each of which is further subdivided into $r$ tiny arcs, one containing $(1-\alpha)s^{2s-i}$ vertices with colour $c_p$, and the other $r-1$ each containing $\alpha s^{2s-i}/(r-1)$ vertices with one of the other $r-1$ colours. For every $q\in[r]$, let the tiny arc in $\lambda_{i,(j-1)r+p,t}$ with colour $c_q$ be $\lambda_{i,(j-1)r+p,t,q}$. Let the majority colour of the medium arc $\ell_{i',j'}$ be $c\in\{c_1,\ldots,c_r\}$, and note that $\ell_{i',j'}$ is also similarly subdivided into tiny arcs like above. 

For every $p,q\in[r]$ and $t\in[s^{s+2i}]$, define $X_{p,q,t}$ as follows.
\[X_{p,q,t}=\begin{dcases}\frac{|\{x\in\lambda_{i,(j-1)r+p,t,q}:x\in V(M'')\}|}{\alpha s^{2s-i'}/(r-1)}, & \text{ if }c_q\not=c,\\\frac{|\{x\in\lambda_{i,(j-1)r+p,t,q}:x\in V(M'')\}|}{(1-\alpha)s^{2s-i'}}, & \text{ if }c_q=c.\end{dcases}\]
Then $\ceil{X_{p,q,t}}$ is the minimum number of tiny arcs in $\ell_{i',j'}$ of colour $c_q$ that contain vertices matched with vertices in $\lambda_{i,(j-1)r+p,t,q}$ by $M''$. Let \[A=\sum_{t=1}^{s^{s+2i}}\sum_{\substack{p,q\in[r]\\c_q\not=c,\;c_p=c_q}}X_{p,q,t},\quad B=\sum_{t=1}^{s^{s+2i}}\sum_{\substack{p,q\in[r]\\c_q=c,\;c_p\not=c_q}}X_{p,q,t},\] \[C=\sum_{t=1}^{s^{s+2i}}\sum_{\substack{p,q\in[r]\\c_q\not=c,\;c_p\not=c_q}}X_{p,q,t},\quad D=\sum_{t=1}^{s^{s+2i}}\sum_{\substack{p,q\in[r]\\c_q=c,\;c_p=c_q}}X_{p,q,t}.\]
Note that by bounding the numerator trivially with $|\lambda_{i,(j-1)r+p,t,q}|$, we get
\[X_{p,q,t}\leq\begin{dcases}\frac{(1-\alpha)s^{2s-i}}{\alpha s^{2s-i'}/(r-1)}=\frac{(r-1)(1-\alpha)}\alpha s^{i'-i}, & \text{ if }c_q\not=c,\;c_p=c_q,\\\frac{\alpha s^{2s-i}/(r-1)}{(1-\alpha)s^{2s-i'}}=\frac{\alpha}{(r-1)(1-\alpha)}s^{i'-i}, & \text{ if }c_q=c,\;c_p\not=c_q,\\\frac{\alpha s^{2s-i}/(r-1)}{\alpha s^{2s-i'}/(r-1)}=s^{i'-i}, & \text{ if }c_q\not=c,\;c_p\not=c_q,\\\frac{(1-\alpha)s^{2s-i}}{(1-\alpha)s^{2s-i'}}=s^{i'-i}, & \text{ if }c_q=c,\;c_p=c_q.\end{dcases}\]
Therefore, $A\leq s^{s+2i}(r-1)^2(1-\alpha)s^{i'-i}/\alpha=(r-1)^2(1-\alpha)s^{s+i+i'}/\alpha$, and similarly $B\leq\alpha s^{s+i+i'}/(1-\alpha)$, $C\leq(r-1)^2s^{s+i+i'}$, and $D\leq s^{s+i+i'}$.

Using the definition of $X_{p,q,t}$, the number of vertices in $\overline{\ell}_{i,j}\cup M'(\overline{\ell}_{i,j})$ that is covered by the matching $M''$ is
\[\frac{2\alpha s^{2s-i'}}{r-1}(A+C)+2(1-\alpha)s^{2s-i'}(B+D).\]
Recall that $|\overline{\ell}_{i,j}|=rs^{3s+i}$. Now, we find a lower bound for $|M'(\overline{\ell}_{i,j})|$. For each $p,q\in[r]$ and $t\in[s^{2s+i}]$, from definition, at least $X_{p,q,t}$ tiny arcs in $\ell_{i',j'}$ of colour $c_q$ contains vertices matched with vertices in $\lambda_{i,(j-1)r+p,t,q}$ by $M''$. Therefore, $M'(\lambda_{i,(j-1)r+p,t,q})$ contains at least $X_{p,q,t}-2$ full small arcs of the form $\lambda_{i',j',t'}$, with $t'\in [s^{s+2i'}]$, so $|M'(\lambda_{i,(j-1)r+p,t,q})|\geq(X_{p,q,t}-1)s^{2s-i'}$. Note that for distinct triples $(p,q,t)$, $M'(\lambda_{i,(j-1)r+p,t,q})$ are disjoint subsets of $M'(\overline{\ell}_{i,j})$, so we have
\[|M'(\overline{\ell}_{i,j})|\geq(A+B+C+D)s^{2s-i'}-2r^2s^{s+2i}s^{2s-i'}=(A+B+C+D)s^{2s-i'}-2r^2s^{3s+2i-i'}.\]
Thus, 
\begin{align*}
\eff(\overline{\ell}_{i,j})&\leq\dfrac{\dfrac{2\alpha s^{2s-i'}}{r-1}(A+C)+2(1-\alpha)s^{2s-i'}(B+D)}{(A+B+C+D)s^{2s-i'}-2r^2s^{3s+2i-i'}+rs^{3s+i}}\\
&=\dfrac{\dfrac{2\alpha}{r-1}(A+C)+2(1-\alpha)(B+D)}{(A+B+C+D)+rs^{s+i+i'}-2r^2s^{s+2i}}=:\eff(A,B,C,D).
\end{align*}

Next, we bound $\eff(A,B,C,D)$ from above. First note that
\[\eff(A,B,C,D)=2(1-\alpha)-\dfrac{\left(2(1-\alpha)-\dfrac{2\alpha}{r-1}\right)(A+C)+2(1-\alpha)\left(rs^{s+i+i'}-2r^2s^{s+2i}\right)}{(A+B+C+D)+rs^{s+i+i'}-2r^2s^{s+2i}}.\]
Recall that $\alpha=1-1/\sqrt r$ and $i'>i$, so $2(1-\alpha)>2\alpha/(r-1)$ and $rs^{s+i+i'}>2r^2s^{s+2i}$ for sufficiently large $s$. Thus, increasing $B$ and $D$ increases $\eff(A,B,C,D)$, and so
\begin{align*}
\eff(A,B,C,D)&\leq\eff\left(A,\dfrac{\alpha s^{s+i+i'}}{1-\alpha},C,s^{s+i+i'}\right)\\
&=\dfrac{\dfrac{2\alpha}{r-1}(A+C)+2s^{s+i+i'}}{(A+C)+\left(r+\dfrac{1}{1-\alpha}\right)s^{s+i+i'}-2r^2s^{s+2i}}\\
&=\dfrac{2\alpha}{r-1}+\dfrac{\left(2-\dfrac{2\alpha}{r-1}\left(r+\dfrac{1}{1-\alpha}\right)\right)s^{s+i+i'}+\dfrac{4\alpha}{r-1}r^2s^{s+2i}}{(A+C)+\left(r+\dfrac{1}{1-\alpha}\right)s^{s+i+i'}-2r^2s^{s+2i}}\\
&=\dfrac{2\alpha}{r-1}+\dfrac{\dfrac{4\alpha}{r-1}r^2s^{s+2i}}{(A+C)+\left(r+\dfrac{1}{1-\alpha}\right)s^{s+i+i'}-2r^2s^{s+2i}},
\end{align*}
where the last equality follows as the coefficient of $s^{s+i+i'}$ in the numerator simplifies to 0 with the choice of $\alpha=1-1/\sqrt r$. 

Finally, using that $s^{s+2i}\leq s^{s+i+i'}/s=o(s^{s+i+i'})$ as $s$ tends to infinity, we have that 
\[\eff(\overline{\ell}_{i,j})\leq\eff(A,B,C,D)\leq\dfrac{2\alpha}{r-1}+o(1)=\dfrac{2}{r+\sqrt r}+o(1).\]
Since this is true for every bonded pair $(i,i')$ and every $j\in[s^{2s-i}]$ that is $(i,i')$-regular, $|V(M'')|\leq\left(\frac2{r+\sqrt r}+o(1)\right)N$. Therefore, the number of edges in $M$, which is a maximum cyclic nested matching in $H$, is at most $|V(M'')|/2+o(N)\leq\left(\frac1{r+\sqrt r}+o(1)\right)N$, as claimed. 
\end{proof}

\section{Concluding Remarks}
Bašić, Damnjanović, Stevanović, and Stošić~\cite{Bavsic2026SomeRO} defined the following much more general notion of permutation Ramsey numbers. 

Let $H_1,\ldots,H_k$ be ordered graphs, and let $\Gamma_i$ be a subgroup of the symmetric group on $V(H_i)$ for every $i\in[k]$. Then, the \textit{permutational Ramsey number} $R(H_1^{\Gamma_1},\ldots,H_k^{\Gamma_k})$ is defined to be the smallest natural number $N$ so that whenever the ordered complete graph $K_N^<$ is $k$-coloured, there exists $i\in[k]$ and $\sigma\in\Gamma_i$, such that there is a colour-$i$ copy of $\sigma(H_i)$, the ordered graph obtained by permuting the vertices of $H_i$ with $\sigma$.

The usual Ramsey numbers, the ordered Ramsey numbers, and the cyclic Ramsey numbers correspond to the cases when each $\Gamma_i$ is the full permutation group, the trivial subgroup, and the subgroup of cyclic shifts, respectively. Bašić et al.~\cite{Bavsic2026SomeRO} also defined the \textit{reflective Ramsey numbers} and the \textit{dihedral Ramsey numbers}, corresponding to the cases when each $\Gamma_i$ is the subgroup generated by the reversal permutation, and the dihedral group, respectively.

During the preparation of this paper, a preprint of Damnjanović and {\DJ}or{\dj}ević~\cite{damnjanovic2026computation} was posted on arXiv. Also using the Kissat SAT solver, the authors systematically computed small reflective and dihedral Ramsey numbers, and formulated six new conjectures. As summarized in Table~\ref{table2} below, five of these conjectures on dihedral Ramsey numbers follow directly from five corresponding conjectures in~\cite{Bavsic2026SomeRO} that we resolved in this paper.

\begin{table}[h!]
\begin{center}\renewcommand{\arraystretch}{1.3}
\resizebox{\textwidth}{!}{%
\begin{tabular}{|c c c c c c|}
\hline
$H_1$ & $H_2$ & Cyclic Conjecture & Dihedral Conjecture & LB Reference & Proof Reference\\
\hline\hline
$\AP_k$ & $\AP_\ell$ & {\cite[Conjecture~4.8]{Bavsic2026SomeRO}} & {\cite[Conjecture~4.4]{damnjanovic2026computation}} & Proposition~\ref{prop:aplower} & Theorem~\ref{thm:cycapap}\\\hline
$S_k$ & $\AP_\ell$ & {\cite[Conjecture~4.17]{Bavsic2026SomeRO}} & {\cite[Conjecture~4.5]{damnjanovic2026computation}} & Proposition~\ref{prop:aplower} & Theorem~\ref{thm:cycapap}\\\hline
$\AP_k$ & $\MP_\ell$ & {\cite[Conjecture~4.9]{Bavsic2026SomeRO}} & {\cite[Conjecture~4.6]{damnjanovic2026computation}} & Lemma~\ref{lemma:monogen} & Theorem~\ref{thm:cycapmp}\\\hline
$\MC_k$ & $\AP_\ell$ & {\cite[Conjecture~4.14]{Bavsic2026SomeRO}} & {\cite[Conjecture~4.8]{damnjanovic2026computation}} & Lemma~\ref{lemma:4.6} & Theorem~\ref{thm:kmcap}\\\hline
$K_k$ & $\AP_\ell$ & {\cite[Conjecture~4.23]{Bavsic2026SomeRO}} & {\cite[Conjecture~4.9]{damnjanovic2026computation}} & Lemma~\ref{lemma:4.6} & Theorem~\ref{thm:kmcap}\\\hline
\end{tabular}
}
\caption{Five conjectures on dihedral Ramsey numbers from~\cite{damnjanovic2026computation} that are confirmed by our resolution of five corresponding conjectures on cyclic Ramsey numbers in~\cite{Bavsic2026SomeRO}.}\label{table2}
\end{center}
\end{table}

This is due to two simple observations. On the one hand, since the dihedral group contains all cyclic shifts, the dihedral Ramsey number is always at most the corresponding cyclic Ramsey number. On the other hand, note that a dihedral copy of an ordered graph $H$ is either a cyclic copy of $H$, or a reversed cyclic copy of $H$. Thus, if a colouring of $K_N^<$ is invariant under vertex reversal, and is a valid lower bound construction for the cyclic Ramsey number, then it is also a valid lower bound construction for the corresponding dihedral Ramsey number. It turns out that the lower bound constructions for the five corresponding cyclic Ramsey number conjectures are indeed invariant under vertex reversal. 

The remaining~{\cite[Conjecture~4.7]{damnjanovic2026computation}} corresponds to~{\cite[Conjecture~4.13]{Bavsic2026SomeRO}}, which we are not able to resolve. Interestingly, the formulas in these two conjectures differ in their ceiling and floor signs, so the resolution of one may not necessarily imply the other. For these two conjectures, as well as~{\cite[Conjecture~4.22]{Bavsic2026SomeRO}}, we do not even have general lower bound constructions that match the conjectured values, unlike for~{\cite[Conjecture~4.32]{Bavsic2026SomeRO}}. It would be very interesting to see if any progress can be made on these outstanding conjectures. Finally, our disproof of~{\cite[Conjecture~4.36]{Bavsic2026SomeRO}} also raises the question of determining the true values of $\rc(\NM_{2k},K_\ell)$, for which we do not have a good guess.

\section*{Statement of AI use}
ChatGPT was only used during the initial literature search, which produced the useful references~\cite{MR2038478,MR4442957,MR3697924,moser1993recent}, and for the final spellcheck and proofreading.
The proofs in this paper, as well as the writing and presentation are the work of the authors, who take full responsibility for the content.

\bibliographystyle{abbrv}  
\bibliography{bibliography}
\end{document}